\documentclass[11pt, 3p, authoryear]{elsarticle}

\usepackage{graphicx} 
\usepackage{multirow}

\usepackage{natbib}
\usepackage{optidef}
\usepackage{algorithm}
\usepackage{algpseudocode}
\usepackage{xcolor}
\usepackage{amsfonts}
\usepackage{subcaption}
\usepackage{array}
\usepackage{hyperref}

\algnewcommand\algorithmicnot{\textbf{not}}
\algnewcommand\algorithmicbreak{\textbf{break}}
\algdef{SE}[IF]{IfNot}{EndIf}[1]{\algorithmicif\ \algorithmicnot\ #1\ \algorithmicthen}{\algorithmicend\ \algorithmicif}%

\title{Acceleration Techniques for Learning Optimal Classification Trees with Integer Programming}
\author[1]{Mitchell Keegan\texorpdfstring{\corref{cor1}}{}}
\ead{mitchk.academic@gmail.com}
\author[2]{Michael Forbes}
\author[1]{Paul Corry}
\author[1]{Mahdi Abolghasemi}

\cortext[cor1]{Corresponding author\\Email address: mitchell.keegan@qut.edu.au}

\affiliation[1]{organization={School of Mathematical Sciences, Queensland University of Technology},
postcode={4000},
city={Brisbane},
country={Australia}}

\affiliation[2]{organization={School of Maths and Physics, University of Queensland},
postcode={4072},
city={Brisbane},
country={Australia}}

\begin{document}

\begin{abstract}
    Decision trees are a popular machine learning model which are traditionally trained by heuristic methods. Massive improvements in computing power and optimisation techniques has led to renewed interest in learning globally optimal decision trees. Empirical evidence shows that optimal classification trees (OCTs) have better out-of-sample performance than heuristic methods. The dominant optimisation paradigms for training OCTs are mixed-integer programming (MIP) and dynamic programming (DP). MIP formulations offer flexibility in the objectives and constraints that are modelled, but suffer from poor scaling in the size of the training dataset and the maximum tree depth. DP models represent the state of the art in scaling for OCTs, but lack some of the flexibility of MIP models. In this paper we present progress on using advanced integer programming methods to integrate ideas from DP models into MIP formulations to begin bridging the scaling gap. Using the existing BendOCT model from the literature as a base model, we introduce valid inequalities, cutting planes, and a primal heuristic to improve the scaling of MIP formulations. We show that these techniques significantly improve the ability of BendOCT to find provably optimal solutions over a wide range of datasets.
\end{abstract}


\maketitle

\section{Introduction}

Decision trees are popular 
machine learning models for classification and regression which work by partitioning the feature space and assigning predictions in each partition. Their tree structure is easy to visualise and as such they are an attractive model in domains like healthcare, finance, and criminal justice in which interpretability is highly valued \citep{rudinStopExplainingBlack2019}. Finding optimal decision trees is $\mathcal{NP}$-hard \citep{hyafilConstructingOptimalBinary1976}, which historically has motivated the construction of decision trees by top down heuristics. The last decade has seen renewed interest in globally optimal trees, beginning with the work of \citet{bertsimasOptimalClassificationTrees2017} who argued that advances in computing power and the increasing sophistication of optimisation techniques merited a re-evaluation of the tractability of learning optimal decision trees. More recently \texorpdfstring{\cite{lindenOptimalGreedyDecision2024}}{[Linden et al., 2024]} found that optimal classification trees (OCTs) are both consistently smaller and improve out-of-sample accuracy by nearly 2\% on average relative to heuristic methods in comprehensive experiments, providing strong evidence for the effectiveness of optimal classification trees. Another advantage over heuristic methods is that the modelling frameworks used to learn optimal trees can incorporate complicated objectives and constraints. Unfortunately, finding optimal trees is highly non-trivial, necessitating complicated optimisation procedures.

The two dominant optimisation paradigms for optimal classification trees are mixed-integer programming (MIP) and dynamic programming\footnote{We use the term dynamic programming loosely to refer to a variety of closely related methodologies.} (DP). A key research focus is scalability with respect to the size of the training dataset and the maximum tree depth. On this front DP represents the state of the art, scaling to datasets with orders of magnitude more samples than MIP formulations. The drawback of DP formulations is their reliance on a strict assumption of the independence of subtrees which restricts their modelling power. The advantage of MIP formulations is their modelling flexibility which allows them to model any linear objective or constraint, including global objectives.  Improving the scaling of MIP formulations is key to exploiting this flexibility.

To date, the primary avenue of research for improving MIP scaling has been developing stronger formulations. A common issue is the use of big-M constraints which are known to create weak relaxations; recent work avoids big-M constraints by specialising to binary data or by Benders decomposition which has improved scaling \citep{aghaeiStrongOptimalClassification2025, boutilierOptimalMultivariateDecision2023}. However, MIP formulations still lag far behind DP. We argue that existing formulation do not leverage the full power of mixed-integer programming and can be improved upon by integrating ideas from state of the art DP methods. One key to the success of DP for learning OCTs is the use of bounds to prune suboptimal regions of the search space. Such bounds are often inferred by reasoning about the structure of optimal trees, a common example being equivalent points. These are samples which are identical in the feature space yet have different classes; they are guaranteed to be routed into the same leaf where they cannot all be simultaneously classified correctly which implies an upper bound on the classification score of those samples. This is commonly exploited to great effect to tighten DP bounds \cite[See][for example]{huOptimalSparseDecision2019}, but has not been explored in the MIP setting even though it could significantly tighten relaxations. The intention of this work is to use advanced integer programming techniques to incorporate such structure and accelerate the convergence of MIP formulations to optimal trees. 

We specialise our acceleration techniques to the BendOCT model introduced by \cite{aghaeiStrongOptimalClassification2025}. It is among the best performing MIP models in the literature and its use of Benders decomposition is convenient for applying the proposed acceleration techniques. Larger models generally have slower relaxation solve times and risk being unable to progress far into the branch and bound tree which can render our methods ineffective. We study OCTs which maximise accuracy subject to a penalty on the number of leaf nodes, with the intention that the proposed techniques can be extended to more general objectives and constraints in future works. We note that conceptually our techniques apply to other MIP formulations, with modifications to suit the decision variables used. Accelerating the convergence of MIP formulations is critical in unlocking their potential for practical applications. The primary contributions of this work are summarised as follows:

\begin{itemize}
    \item A derivation of BendOCT by logic-based Benders decomposition (LBBD) which differs from the original derivation of \cite{aghaeiStrongOptimalClassification2025}. The intuition underlying the derivation is extended to provide a novel strengthening of the Benders cuts.
    \item A solution polishing primal heuristic, a class of valid inequalities which exploit equivalent points, and a class of cutting planes which exploit a subroutine for quickly finding optimal depth two subtrees.
    \item Computational experiments across 33 datasets validating the effectiveness of the proposed techniques. Compared to the original BendOCT with a one hour time limit our acceleration techniques allow us to solve the same amount of instances to optimality within 12 seconds, and over twice as many within the one hour time limit.
\end{itemize}

The remainder of the paper is organised as follows. Section \ref{sect:Literature Review} reviews related literature. Section \ref{sect:Problem Setting} introduces notation to model OCTs and formally states the problem being studied. Section \ref{sect:BendOCT} derives the BendOCT model by LBBD. Section \ref{sect:Acceleration Techniques} introduces a strengthening of the Benders cuts, a primal heuristic, a class of valid inequalities, and a class of cutting planes which accelerate the convergence of BendOCT. Section \ref{sect: Computational Experiments} describes our experimental setup and presents the results of computational experiments which validate the effectiveness of the proposed acceleration techniques across 33 datasets. Section \ref{sect:Conclusion} concludes the paper. Extended results are available in \ref{Appendix: Extended Results}.

\section{Related Work}\label{sect:Literature Review}
The most popular optimisation frameworks for learning optimal classification trees are MIP and DP. The types of classification trees they learn and the assumptions they make broadly vary along three dimensions. The first is whether they operate on binary, categorical, continuous features, or some combination of the three. Those specialised to binary data rely on discretisation to handle categorical and continuous features at the cost of optimality with respect to the original features. The second is the form of regularisation. This is done in various ways, the most common being a limit on the tree depth (called balanced trees if grown to said limit), a hard constraint on the number of leaf nodes, or by a penalty on the number of leaf nodes. The third is whether branch node splitting rules are univariate or multivariate. Multivariate splits are more expressive, resulting in compact trees which can rival the performance of tree ensembles, but are significantly harder to optimise and sacrifice much of the interpretability of decision trees.

Finding optimal classification trees is $\mathcal{NP}$-hard \citep{hyafilConstructingOptimalBinary1976}, which historically has motivated using greedy heuristics for constructing optimal trees. The most popular of these algorithms build trees in a top down approach, choosing splitting rules at each branch node by locally maximising some metric of the class separation as a proxy for accuracy. The most prominent examples of this approach are CART \citep{breimanClassificationRegressionTrees2017}, ID3 \citep{quinlanInductionDecisionTrees1986}, and C4.5 \citep{quinlanC45ProgramsMachine1993}. 

Recent interest in optimal trees begins with the OCT-BD\footnote{The authors call their model OCT, we only name it OCT-BD here for convenience} MIP formulation of \cite{bertsimasOptimalClassificationTrees2017}. They model univariate and multivariate splits on continuous features and allow regularisation by penalising the
number of branch nodes. Their model showed promise in learning small trees but suffers from poor scaling which, despite significant progress since, remains a weak point of MIP formulations. \citet{verwerLearningOptimalClassification2019} proposed BinOCT for learning balanced trees on binary data. Instead of optimising over binary variables corresponding to feature thresholds BinOCT optimises over binary encodings of the thresholds. \citet{gunlukOptimalDecisionTrees2021} developed a formulation specialised for categorical data which optimises fixed tree topologies. They develop valid inequalities to strengthen the relaxation of the formulation, although these are implied by the model constraints, not inferred from the problem structure.

One source of the poor scaling of these formulations is the inability of the MIP framework to exploit the independence of subtrees. This represents a fundamental drawback, but conversely allows for greater modelling flexibility. A more addressable source of poor scaling is weak relaxations as a result of big-M constraints. To address this \cite{aghaeiStrongOptimalClassification2025} developed FlowOCT to learn trees with univariate splits on binary data which allows them to avoid using big-M constraints. In this way they develop a stronger formulation while retaining the flexibility of the MIP framework. They further show that FlowOCT is amenable to Benders decomposition and extend it to BendOCT. In a similar vein \cite{boutilierOptimalMultivariateDecision2023} proposed S-OCT for learning trees with multivariate splits on continuous features. They avoid big-M constraints by a Benders decomposition in which the master problem constructs a routing of samples through the tree and feasibility subproblems ensure that there exists a multivariate splitting rule which can realise the proposed sample routing. Information is returned from infeasible subproblems by shattering inequalities which encode subsets of samples which cannot be separated by any multivariate split. The shattering inequalities are valid inequalities for MIP formulations with multivariate splits, an idea which is expanded upon by \citet{michiniPolyhedralStudyMultivariate2024} in a polyhedral study of S-OCT.

FlowOCT/BendOCT and S-OCT are the most effective formulations for univariate and multivariate splits respectively. OCT-BD also remains popular for a variety of extensions which leverage the flexibility of mixed-integer programming, most notably prescriptive trees \citep{bertsimasOptimalPrescriptiveTrees2019, ikedaPrescriptivePriceOptimization2023}, robust classification \citep{vosRobustOptimalClassification2022}, learning policies for Markov Decision Processes \citep{vosOptimalDecisionTree2023}, and smart predict-then-optimise \citep{elmachtoubDecisionTreesDecisionMaking2020}. Similarly FlowOCT has been extended to consider prescriptive trees \citep{joLearningOptimalPrescriptive2023}, learning fair decision trees \citep{joLearningOptimalFair2023}, and robust classification \citep{justinLearningOptimalClassification2023}. While effective these extensions would benefit from better scaling to allow application to larger datasets and deeper trees.

The second and more popular method for learning optimal classification trees is DP. DP formulations scale extraordinarily well owing to linear time complexity in the number of training samples. This allows depth four optimal trees to be found in minutes for datasets with up to a million samples and 150 features \citep{lindenOptimalGreedyDecision2024}. MIP formulations by comparison only scale to thousands of samples. Most closely related to our work is a series of models beginning with OSDT \citep{huOptimalSparseDecision2019}, its generalisation to consider imbalanced data and rank statistics \citep{linGeneralizedScalableOptimal2020}, and using reference tree ensembles to guide the exploration of the search space \citep{mctavishFastSparseDecision2022a}. Key to the success of OSDT and its successors is the use of analytic bounds to aggressively prune the search space. In this work we adapt one such bound based on equivalent points into the MIP setting but believe that this avenue holds further opportunity. In a similar vein \citet{aglinLearningOptimalDecision2020} introduced DL8.5 which caches previously seen subtrees to prune similar subtrees without fully exploring them. \cite{demirovicMurTreeOptimalDecision2022} develop bounds which prune the search space by considering similar subtrees which have been solved to optimality and introduce a specialised subroutine for finding optimal depth two subtrees which significantly accelerates convergence. Other recent advancements include Branches \citep{chaoukiBranchesFastDynamic2024}, which employs a novel purification bound to prune the search space. \citet{demirovicBlossomAnytimeAlgorithm2023} improve on the search strategy of DL8.5. \citet{vanderlindenNecessarySufficientConditions2023} introduce a generalised DP framework for learning optimal trees. A drawback of these methods is that bar a few examples \citep{chaoukiBranchesFastDynamic2024, britaOptimalClassificationTrees2025, mazumderQuantBnBScalableBranchandBound2022}, DP formulations only operate on binary features which necessitates a discretisation of continuous features thereby sacrificing optimality. In all cases they cannot learn trees with multivariate splitting rules. They are also limited in the objectives and constraints which can be modelled by the independence of subtrees, we refer to \cite{vanderlindenNecessarySufficientConditions2023} for a formal treatment of the limits of DP formulations. 

Outside of the aforementioned there has been little research on incorporating valid inequalities or cutting planes into MIP formulations, and to the knowledge of the authors none that explicitly exploit the structure of optimal decision trees. DP formulations explicitly leverage this structure to prune the search space and it is fundamental to their success. A partial reason for this may be that the underlying logic is natural to incorporate into DP formulations but can be awkward to model in a MIP setting. Nevertheless, it is a fruitful avenue for improving the scaling of MIP formulations.

\section{Problem Setting}\label{sect:Problem Setting}
In this section we introduce notation for modelling decision trees. We then state a general form of the problem of learning optimal classification trees and formally describe the objective of this work.

\subsection{Notation}
\begin{figure}[ht]
    \centering
    \includegraphics[width=0.5\linewidth]{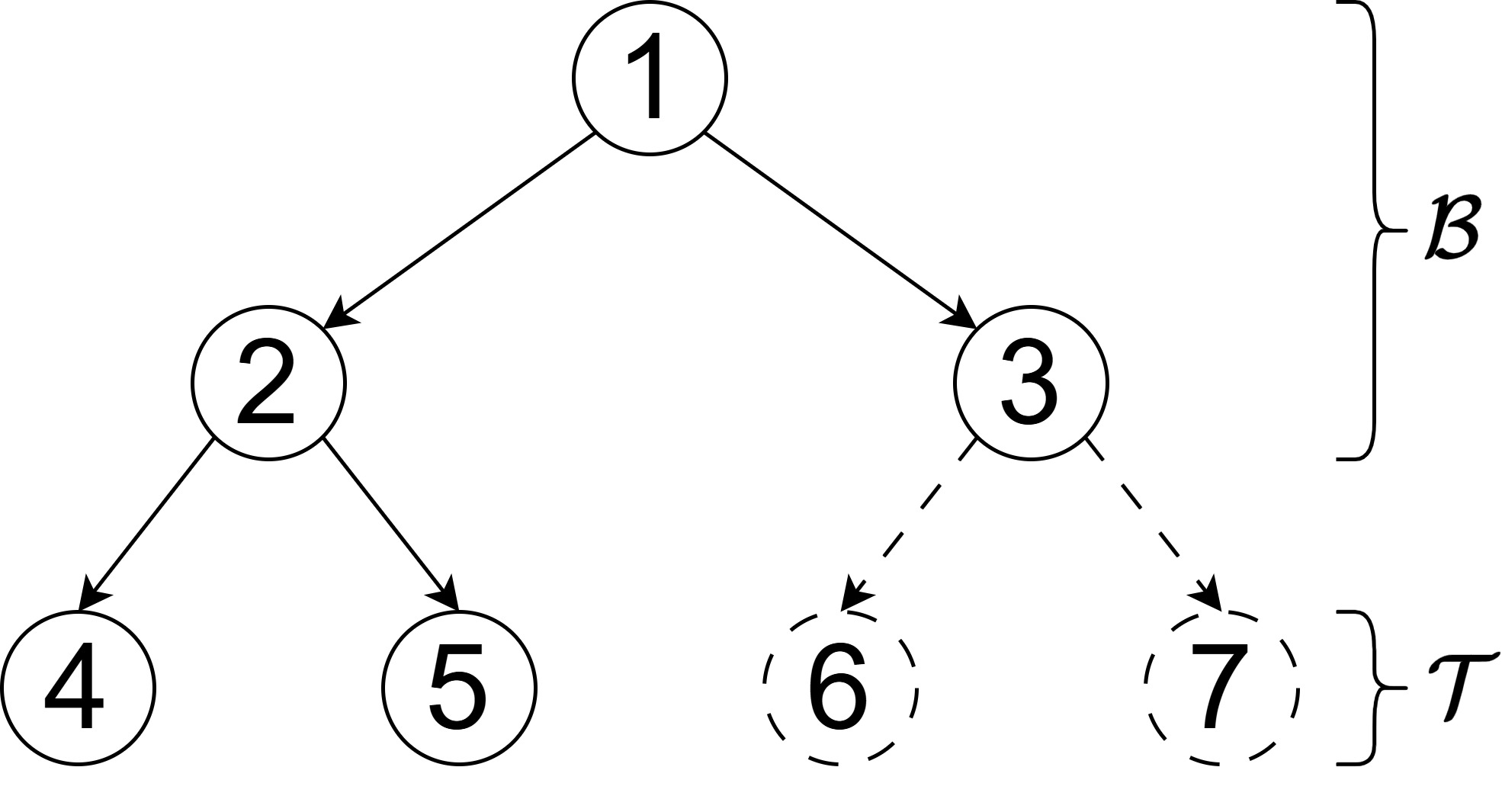}
    \caption{Basic decision tree notation. The available nodes are partitioned into internal nodes $\mathcal{B}$ and terminal nodes $\mathcal{T}$. Nodes 1 and 2 have been designated as branch nodes while nodes 3, 4, and 5 are designated as leaf nodes. Nodes 6 and 7 are cut off.}
    \label{fig:Tree Notation}
\end{figure}

This work considers binary decision trees with maximum depth $D$. Nodes are ordered in breadth-first search ordering as in Figure \ref{fig:Tree Notation} which displays a depth two tree. Note that the tree is sparse in the sense that not all available nodes are used. The set of available nodes is partitioned into internal nodes $\mathcal{B} := \{1,\hdots,2^D-1 \}$ and terminal nodes $\mathcal{T} := \{2^D,\hdots,2^{D+1}-1\}$. In a decision tree each available node can have one of three designations:
\begin{itemize}
    \item Branch node - A node with two child nodes and a binary branching rule.
    \item Leaf node - A node at which a prediction is made.
    \item Cut node - Any nodes which are descendant from a leaf node and are neither a branch node nor a leaf node.
\end{itemize}

Internal nodes in $\mathcal{B}$ can take any of the three designations, except for the root node which must be either a branch node or a leaf node. Terminal nodes in $\mathcal{T}$ must be a leaf node if none of their ancestors are a leaf node, otherwise they are cut. 

For any internal node $n \in \mathcal{B}$, we take the left and right children to be $l(n) = 2n$ and $r(n) = 2n+1$ respectively. The parent of node $n$ is $a(n) = \lfloor \frac{n}{2} \rfloor$ except for the root node which does not have a parent. Denote the set of ancestors of $n$, not including $n$, as $\mathcal{A}(n)$. Further let $\mathcal{A}(n)$ be partitioned by $\mathcal{A}^L(n)$ and $\mathcal{A}^R(n)$ which are the ancestors of $n$ which branch left and right respectively on the path to $n$. For example, in Figure \ref{fig:Tree Notation} $\mathcal{A}^L(5) = \{1\}$, $\mathcal{A}^R(5) = \{2\}$, and $\mathcal{A}(5) = \{1,2\}$. Similarly, the descendants of node $n$, not inclusive of $n$, are denoted $\mathcal{D}(n)$. For convenience we define the descendants and ancestors of node $n$, inclusive of $n$, as $\mathcal{D}^+(n)$ and $\mathcal{A}^+(n)$ respectively. For two nodes $n_1$ and $n_2$ such that $n_2 \in \mathcal{D}^+(n_1)$ we define the distance between them, $dist(n_1,n_2)$, to be the number of edges traversed on a path between the two nodes. The height of a node $n$, denoted $height(n) = \max_{n_d \in \mathcal{D}^+(n)} dist(n,n_d)$, is its distance from its furthest descendant leaf node\footnote{The height may be defined on either the full tree of available nodes, in which case the descendant leaf nodes are the terminal nodes $\mathcal{T}$, or a sparse tree where the descendant leaf nodes are those assigned as such. When not clear from context we will explicitly distinguish between the two.}. Similarly, the depth of a node $n$ is $dist(1,n)$.

In this work decision trees are trained on and make predictions on datasets with binary features and discrete classes. For datasets with continuous or categorical features a binary encoding of the data is used as described in Section \ref{Binary Data Encoding}. The set of indices of samples is denoted by $\mathcal{I}$, the set of binary features by $F$, and the set of target classes by $\mathcal{K}$. The feature data is denoted $X \in \{0,1\}^{|\mathcal{I}| \times |F|}$ where $x_f^i \in \{0,1\}$ is the value of feature $f \in F$ in sample $i \in \mathcal{I}$. Similarly, the target vector is $Y \in \mathcal{K}^{|\mathcal{I}|}$ where $y^i \in \mathcal{K}$ is the class of sample $i$. If a node branches on feature $f$, then any data sample such that $x_f^i = 0$ ($x_f^i = 1$) is routed to the left (right) child node. A sample $i$ is routed through a decision tree according to the branch rules as described and terminates in some leaf node. It is correctly classified if the leaf node predicts class $y^i$ and misclassified otherwise. A feasible assignment of available nodes is one in which the following hold:
\begin{itemize}
    \item Each internal node must either be a branch node, a leaf node, or one of its ancestor must be a leaf node.
    \item Exactly one node on the paths from the root node to each of the terminal nodes must be a leaf node.
\end{itemize}
A decision tree consists of a feasible assignment of all available nodes $n \in \mathcal{B} \cup \mathcal{T}$ as branch, leaf, or cut nodes, a single feature $f \in F$ representing a branching rule for each assigned branch node, and a prediction $k \in \mathcal{K}$ for each assigned leaf node.

\subsection{Problem Statement}
Given some decision tree $\mathcal{DT}$, denote the predicted class of the $i^{th}$ sample as $\mathcal{DT}(x^i)$, the number of leaves in $\mathcal{DT}$ as $H(\mathcal{DT})$, and the set of feasible decision trees as $\mathcal{DT}^F$. Formally we solve the optimisation problem of finding a feasible decision tree which maximises the accuracy over training data $(X,y)$ subject to a penalty $\lambda$ on the number of leaf nodes in the tree:
\begin{maxi}[s]
    {\mathcal{DT} \in \mathcal{DT}^F}
    {\frac{|\{i \in \mathcal{I} : \mathcal{DT}(x^i) = y^i \}|}{|\mathcal{I}|} - \lambda H(\mathcal{DT})}{\label{model:General OCT Problem}}
    {}
\end{maxi}

Many MIP formulations of (\ref{model:General OCT Problem}) and its variants have been proposed. Instead of proposing a new formulation we aim to develop techniques which can accelerate the convergence of existing formulations. We build our methods on the BendOCT model, but conceptually our techniques should apply to other MIP formulations, particularly those with univariate decision rules, with modifications to suit the decision variables used. We focus our efforts on classification trees which maximise accuracy, but note that MIP formulations generalise well to any metric which can be expressed in linear terms. MIP frameworks are also flexible in the constraints placed on the decision tree structure $\mathcal{DT}^F$, for example limiting the number of unique features used or disallowing specific feature combinations. $\mathcal{DT}^F$ can also be made data dependent, one such example being a constraint on the minimum number of samples in leaf nodes which is a common regularisation method for decision trees. An interesting direction for future work is to extend the acceleration techniques introduced in Section \ref{sect:Acceleration Techniques} to more general objectives and constraints.

The regularisation parameter $\lambda$ can be interpreted as the required marginal improvement in accuracy for each additional leaf node. Some form of tunable regularisation parameter which controls the size, and by extension the complexity, of the decision tree is required to prevent overfitting and improve generalisation performance. The maximum depth $D$ acts as a tunable parameter but in practice is too coarse for usage in isolation. An alternative is to place a hard limit on the number of leaf nodes in the tree by adding constraint $H(\mathcal{DT}) \leq C$ for some $C$. As per \citet{lombardiAnalysisRegularizedApproaches2021} regularisation by a hard constraint is more expressive than penalising the objective in the sense that for any $\lambda \in \mathbb{R}^{\geq 0}$ there exists a $C \in \{ 1,\hdots ,2^D \}$ which yields the same optimal objective, while the converse does not hold. We choose regularisation by objective penalisation because it is often used in related works and some of the acceleration techniques we introduce do not work with a hard constraint on tree size.


\section{The BendOCT MIP Formulation}\label{sect:BendOCT}
BendOCT is a MIP model for constructing optimal classification trees introduced by \citet{aghaeiStrongOptimalClassification2025}. They first develop a flow based formulation named FlowOCT which models the path of each data sample through the decision tree by flow variables and links the structure of tree to the sample flow such that only correctly classified samples have a non-zero flow. BendOCT is then developed by noting that after fixing variables relating to the tree structure the optimisation problem decomposes into $|\mathcal{I}|$ independent linear programming (LP) subproblems. While the resulting LP subproblems can be solved by a generic solver and Benders cuts generated by standard duality theory, the authors instead note that each subproblem constitutes a max-flow problem. This allows them to exploit max-flow/min-cut duality to generate Benders cuts by an efficient minimum cut routine tailored to the structure of the subproblem flow graphs.

In this work we propose an alternative derivation of BendOCT by logic-based Benders decomposition (LBBD). Introduced by \citet{hookerLogicbasedBendersDecomposition2003}, LBBD is a generalisation of Benders decomposition in which the subproblems take a more general form than a linear program and Benders cuts are derived by logical inference on the subproblems instead of by solving the dual of the subproblem LP. We begin by defining decision variables which describe the structure of the decision tree:
\begin{itemize}
    \item $b_{nf} \in \{0,1\}$ is equal to one if internal node $n \in \mathcal{B}$ is a branch node and branches on feature $f \in F$. 
    \item $p_n \in \{0, 1\}$ is equal to one if node $n \in \mathcal{B} \cup \mathcal{T}$ is a leaf node.
    \item $w_k^n \in \{0,1\}$ is equal to one if node $n \in \mathcal{B} \cup \mathcal{T}$ predicts class $k \in \mathcal{K}$.
\end{itemize}

A valid decision tree is one which obeys the following constraints:

\begin{itemize}
    \item Each internal node must either branch on a feature $f \in F$, be a leaf node predicting a class $k \in \mathcal{K}$, or one of its ancestors must be a leaf node.
    \item A prediction must be made at exactly one node on the paths from the root node to each of the terminal nodes.
\end{itemize}

Any nodes which are cut off will have all associated decision variables set to zero. Outside of enforcing a valid tree structure, any formulation must link the structure of the tree to the routing and classification of samples. To this end, denote $\mathcal{DT}(b,p,w)$ to be the tree defined by decision variables $b$, $p$, and $w$, then define:

\begin{equation}\label{BendOCT Theta}
    \Theta_i(b,p,w) = \begin{cases}
        1,\; \text{if sample $i$ is classified correctly by tree } \mathcal{DT}(b,p,w) \\
        0,\; \text{otherwise}
    \end{cases}
\end{equation}

$\Theta_i(b,p,w)$ is the classification score of the $i^{th}$ sample and can be evaluated by following the path of sample $i$ through tree $\mathcal{DT}(b,p,w)$ and checking if it is correctly classified in a leaf node. A minimal model is then:

\begin{maxi!}[s]
    {b,p,w}
    {\frac{1}{|\mathcal{I}|}\sum_{i \in \mathcal{I}} \Theta_i(b,p,w) - \lambda \sum_{n \in \mathcal{B} \cup \mathcal{T}} p_n}
    {\label{model:BendOCT Minimal}}
    {}
    \addConstraint{\sum_{f \in \mathcal{F}} b_{nf} + p_n + \sum_{n_a \in \mathcal{A}(n)} p_{n_a}}{=1,\quad}{\forall n \in \mathcal{B} \label{const:BendOCTMinimal branch variables}}
    \addConstraint{\sum_{k \in \mathcal{K}} w_k^n}{=p_n, \quad}{\forall n \in \mathcal{B} \cup \mathcal{T} \label{const:BendOCTMinimal prediction variables}}
    \addConstraint{p_n + \sum_{n_a \in \mathcal{A}(n)} p_{n_a}}{=1,\quad}{\forall n \in \mathcal{T} \label{const:BendOCTMinimal path predictions}}
    \addConstraint{b_{nf}}{\in \{0,1\} \quad}{\forall n \in \mathcal{B}, \forall f \in F}
    \addConstraint{p_n}{\in \{0,1\} \quad}{\forall n \in \mathcal{B} \cup \mathcal{T}}
    \addConstraint{w_k^n}{\in \{0,1\} \quad}{\forall n \in \mathcal{B} \cup \mathcal{T}, \forall k \in K \label{const:BendOCTMinimal Last Constraint}}
\end{maxi!}

The objective maximises accuracy subject to a penalty on the number of leaf nodes used. Constraints (\ref{const:BendOCTMinimal branch variables}) enforces that each internal node is a branch node branching on exactly one feature, is a leaf node itself, or has an ancestor which is a leaf node. Constraints (\ref{const:BendOCTMinimal prediction variables}) enforces that each leaf node predicts exactly one class. Constraints (\ref{const:BendOCTMinimal path predictions}) states that on the path from the root node to each terminal node, exactly one node must be designated as a leaf node. This model is clearly not solvable as a MIP owing to the nonlinearity of $\Theta_i(b,p,w)$. As such, new variables $\theta_i$ for $i \in \mathcal{I}$ are introduced to approximate the contribution of $\Theta_i(b,p,w)$ to the objective. Initially these will overapproximate $\Theta_i(b,p,w)$; Benders cuts are then added to refine the approximation. The master problem (MP) is then defined as:

\begin{maxi!}[s]
    {b,p,w,
    \theta}
    {\frac{1}{|\mathcal{I}|}\sum_{i \in \mathcal{I}} \theta_i - \lambda \sum_{n \in \mathcal{B} \cup \mathcal{T}} p_n}
    {\label{model:BendOCT}}
    {}
    \addConstraint{(\ref{const:BendOCTMinimal branch variables}) - (\ref{const:BendOCTMinimal Last Constraint})}{}{}
    \addConstraint{\theta_i}{\in Bend, \quad}{\forall i \in \mathcal{I} \label{const:BendOCT Benders Cuts}}
    \addConstraint{\theta_i}{\in [0,1]}{\forall i \in \mathcal{I}}
\end{maxi!}

$\theta$ is bounded from above so that the model is not unbounded before any Benders cuts are added. Constraints (\ref{const:BendOCT Benders Cuts}) are the Benders cuts inferred from the Benders subproblems which refine the approximations of $\Theta_i(b,p,w)$. At each iteration solving the MP returns an integral solution $(b,p,w,\theta)$. For each sample $i \in \mathcal{I}$, $\Theta_i(b,p,w)$ is evaluated and compared to $\theta_i$ to assess the accuracy of the approximation. If $\theta_i > \Theta_i(b,p,w)$ then the MP has overestimated the classification score of the $i^{th}$ sample and Benders cuts are added to refine the approximation. An important point is that we seek to design Benders cuts which refine the approximation of $\theta_i$ not just for the current tree structure, but also similar tree structures.

\begin{figure}[htb!]
    \centering
    \includegraphics[width=0.6\linewidth]{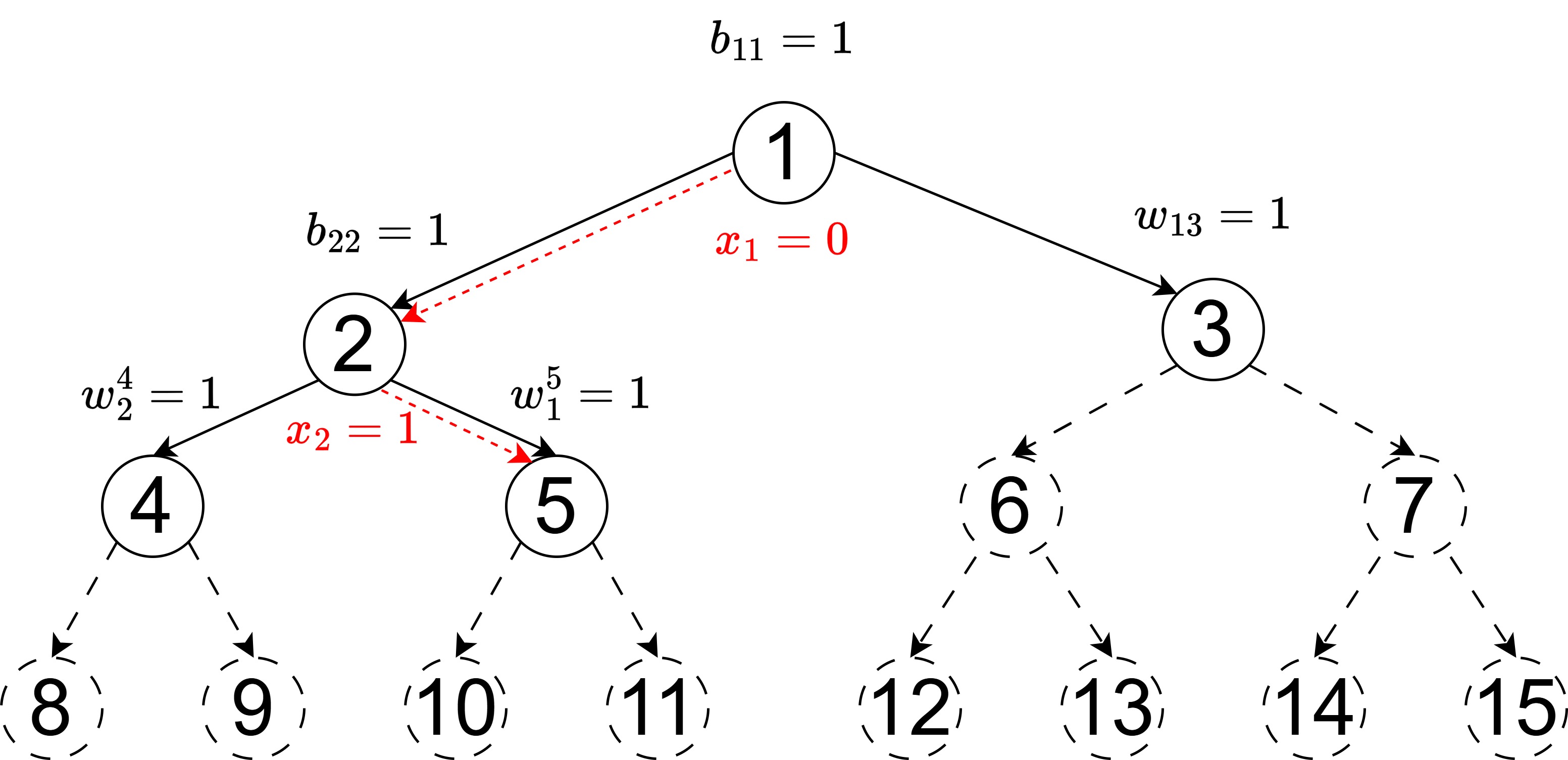}
    \caption{Logic-based Benders cut example tree structure. The sample is routed into node $5$ where it is misclassified.}
    \label{fig:LBBD Example}
\end{figure}

We motivate the structure of the Benders cuts by way of example. Consider a problem with a depth $3$ tree. Take $\mathcal{F} = \{1,2,3,4\}$, $\mathcal{K} = \{1,2\}$, $\mathcal{B} = \{1,2,3,4,5,6,7\}$, and $\mathcal{T} = \{8,9,10,11,12,13,14,15\}$. Assume the following decision variables are equal to one in the MP solution: $b_{11}, b_{22}, w_1^3, w_2^4, w_1^5$. Figure \ref{fig:LBBD Example} shows the resulting tree structure. Consider some sample with $x = (0,1,1,0)$ and $y=2$. It will fall into node $5$ and be misclassified. If $\theta=1$ for that sample, the MP has over estimated $\Theta(b,p,w)$ and a cut should be added to the MP to cut off the solution.

There is an obvious no-good cut, in that $\theta = 1$ is invalid for the given tree structure within the nodes of the tree which the sample passes through, i.e. $b_{11} \land b_{22} \land w_1^5 \implies \theta = 0$. This corresponds to the constraint:

\begin{equation}\label{BendOCT no-good cut}
    \theta \leq (1-b_{11}) + (1 - b_{22}) + (1 - w_1^5)
\end{equation}

This can be immediately strengthened by noting that other branching decisions would also route $x$ into the same leaf node. In this case branching on the $4^{th}$ feature at node $1$ or the $3^{rd}$ feature at node 2 has the same effect, so the no-good cut (\ref{BendOCT no-good cut}) can be strengthened to:

\begin{equation}
    \theta \leq (1-b_{11}-b_{14}) + (1 - b_{22} - b_{23}) + (1 - w_1^5)
\end{equation}

The cut can be further strengthened by noting that $\theta$ should also be forced to zero if nodes $1$ or $2$ predict class $1$:

\begin{equation}
    \theta \leq (1-b_{11}-b_{14} - w_1^1) + (1 - b_{22} - b_{23} - w_1^2 - p_1) + (1 - w_1^5 - p_1 - p_2)
\end{equation}

The additional $p$ variables are required to maintain the logic of the cut. If node $2$ predicts class $1$ then $p_2$ must be subtracted so that the terms associated with node $5$ still sum to zero. Similarly, $2p_1$ is subtracted to account for predictions made at node $1$. By leveraging constraints (\ref{const:BendOCTMinimal branch variables}) and (\ref{const:BendOCTMinimal prediction variables}) it can be shown that this is equivalent to:

\begin{equation}
    \theta \leq (b_{12} + b_{13} + w_2^1) + (b_{21} + b_{24} + w_2^2) + (b_{41} + b_{42} + b_{43} + b_{44} +w_2^5)
\end{equation}

Intuitively this states that $\theta$ would be unbounded if the branch rules route the sample into a different leaf node, or if any of the nodes $1$, $2$, or $5$ predict class $2$. Note that $\theta$ being unbounded for some MP solution only implies that the sample is not strictly misclassified, it does not mean that the sample will be classified correctly. It may be the case that other constraints or Benders cuts simultaneously force $\theta$ to zero. In a feasible solution only $\theta_i$'s which are not forced to zero by any Benders cuts are equal to $1$.

In general, consider some sample $x^i$ which is sorted into leaf $n_l \in \mathcal{B} \cup \mathcal{T}$, where $w_{y^i}^{n_l} = 0$ and $\theta_i = 1$. i.e. $\theta_i$ overestimates $\Theta_i(b,p,w)$. The desired cut will leave $\theta_i$ unbounded if a branch decision is made at any ancestor node $n_a \in \mathcal{A}(n_l)$ or in node $n_l$ itself that would route the sample to a different leaf node, or if any nodes on the path from the root node to $n_l$ are leaf nodes which predict class $y^i$. The following cut encodes this logic:
\begin{equation}\label{Benders Cuts}
    \theta_i \leq \sum_{\substack{n_a \in \mathcal{A}^L(n_l) \\ f \in F : x_f^i = 1}} b_{n_af} + \sum_{\substack{n_a \in \mathcal{A}^R(n_l) \\ f \in F : x_f^i = 0}} b_{n_af} + \sum_{f \in F} b_{n_lf} + \sum_{n_a \in \mathcal{A}^+(n_l)} w_{y^i}^{n_a}
\end{equation}
when $n_l \in \mathcal{B}$ and 
\begin{equation}\label{Benders Cuts Special Case}
    \theta_i \leq \sum_{\substack{n_a \in \mathcal{A}^L(n_l) \\ f \in F : x_f^i = 1}} b_{n_af} + \sum_{\substack{n_a \in \mathcal{A}^R(n_l) \\ f \in F : x_f^i = 0}} b_{n_af} + \sum_{n_a \in \mathcal{A}^+(n_l)} w_{y^i}^{n_a}
\end{equation}
in the special case when $n_l \in \mathcal{T}$. These are exactly the Benders cuts derived by \cite{aghaeiStrongOptimalClassification2025}. A minor improvement is to relax the integrality of the prediction variables (\ref{const:BendOCTMinimal Last Constraint}). After fixing binary solutions for $b$ and $p$, the MP (\ref{model:BendOCT}) decomposes per leaf node with the prediction variables in each leaf node being set to maximise the classification score in the leaf as per the Benders cuts. If there is a unique majority class then the prediction variables which maximise the classification score are obviously binary. If there is not a unique majority class then any solution in which the prediction variables for the majority classes sum to one is optimal, one of these optimal solutions being a binary prediction of one of the majority classes. We implement the Benders cuts as lazy cuts in a callback function.

\section{Acceleration Techniques}\label{sect:Acceleration Techniques}

In this section we detail our core contribution, a number of acceleration techniques which significantly reduce the time taken to find optimal classification trees with the BendOCT model and shrink the optimality gap in instances where optimal solution cannot be found. In Section \ref{sect:Strengthened Benders Cuts} we describe a strengthened version of the Benders cuts which can be applied whenever the sample in the subproblem is routed into a terminal node. Section \ref{sect:Solution Polishing Primal Heuristic} introduces a variation of a subroutine from the dynamic programming literature for finding optimal depth two subtrees which we denote D2S and employs it in a primal heuristic which fixes the upper tree layers from incumbent solutions and optimises the tails of the tree. Section \ref{sect: Equivalent Point Valid Inequalities} introduces a class of valid inequalities which exploits sets of sample which are equivalent or nearly equivalent in the feature space but have differing classes. Section \ref{sect:Path Bound Cutting Planes} introduces cutting planes which exploit the optimality of the D2S subroutine to aggressively prune the search space. For some acceleration techniques we introduce multiple variants which are experimentally evaluated in \ref{Appendix: Extended Results}. 

\subsection{Strengthened Benders Cuts}\label{sect:Strengthened Benders Cuts}
An advantage of deriving the Benders cuts by LBBD is that the intuitive perspective affords a strengthening of the cuts. Assume the $i^{th}$ sample is routed into terminal node $n_l \in \mathcal{T}$. We denote $n_p = a(n_l)$ as the parent node of $n_l$ and assume that $n_l$ is the left child of $n_p$ with sibling node $n_s = r(n_p)$. The strengthening trivially extends to the case where $n_l$ is the right child. Then we can alternatively express the Benders cuts (\ref{Benders Cuts Special Case}) as:

\begin{equation}
    \theta_i \leq \sum_{\substack{n_a \in \mathcal{A}^L(n_p) \\ f \in F : x_f^i = 1}} b_{n_af} + \sum_{\substack{n_a \in \mathcal{A}^R(n_p) \\ f \in F : x_f^i = 0}} b_{n_af} + \sum_{f \in F : x_f^i = 1} b_{n_pf} + \sum_{n_a \in \mathcal{A}^+(n_l)} w_{y^i}^{n_a}
\end{equation}

by separating out the branch variables relating to the parent node $n_p$. The contribution of the sum over $b_{n_pf}$ is to assert that $\theta_i$ can be unbounded if sample $i$ is routed into the sibling node $n_s$ instead of $n_l$ where it was misclassified. We can tighten this condition to assert that $\theta_i$ can be fully unbounded only if the sibling node $n_s$ also correctly predicts class $y^i$. Whenever node $n_s$ misclassifies sample $i$ (i.e. $w_{y^i}^{n_s} = 0$) we can strengthen the cut to:

\begin{equation}
    \theta_i \leq \sum_{\substack{n_a \in \mathcal{A}^L(n_p) \\ f \in F : x_f^i = 1}} b_{n_af} + \sum_{\substack{n_a \in \mathcal{A}^R(n_p) \\ f \in F : x_f^i = 0}} b_{n_af} + \frac{1}{2}\left( \sum_{f \in F : x_f^i = 1} b_{n_pf} + w_{y^i}^{n_s} \right) + \sum_{n_a \in \mathcal{A}^+(n_l)} w_{y^i}^{n_a}
\end{equation}

The logic underlying the strengthening can be extended to other cases. For example if $height(n_l) = 1$, i.e. $n_l$'s children are terminal nodes, we can tighten the condition on branching at $n_l$ instead of making a prediction to branch and making a correct prediction at one of the child nodes. Early experiments on these variants showed no discernible benefit. We note that under the strengthened cuts we no longer relax the prediction variables $w$, since it is no longer true that after fixing binary $(b,p)$ that $w$ will be binary in the relaxation solution. There exist "fractional" predictions which better classify the samples than any binary predictions.

\subsection{Solution Polishing Primal Heuristic}\label{sect:Solution Polishing Primal Heuristic}

Given access to a computationally cheap algorithm which finds strong solutions for small trees, we can develop a primal heuristic which fixes the top layers of a proposed tree and optimises the tails of the tree to attempt to find improved incumbent solutions. Such an algorithm was proposed by \citet{demirovicMurTreeOptimalDecision2022} for depth two classification trees. The subroutine returns an optimal depth two subtree and has $\mathcal{O} (|\mathcal{I}|\cdot |F|^2)$ time complexity and $\mathcal{O} (|F|^2)$ space complexity. We modify the subroutine to return the optimal subtree with a maximum depth of two given regularisation parameter $\lambda$. A complete description is provided in \ref{Appendix:Optimal Depth 2 Tree Subroutine}. For our purposes we assume access to subroutine $D2S(\bar{\mathcal{I}},\lambda)$ which takes as an input a subset of the samples $\bar{\mathcal{I}} \subseteq \mathcal{I}$ and regularisation parameter $\lambda$ and returns the following information about an optimal subtree solution:
\begin{itemize}
    \item $\hat{\mathcal{I}}^+$ - The set of samples in $\bar{\mathcal{I}}$ which are correctly classified.
    \item $\hat{\mathcal{I}}^-$ - The set of samples in $\bar{\mathcal{I}}$ which are incorrectly classified.
    \item $\mathcal{B}_{sub}$ - A set of between zero and three tuples $(n,f)$ where $n$ is a branch node in the subtree and $f$ is the branch feature at node $n$ in the optimal subtree.
    \item $\mathcal{L}_{sub}$ - A set of between one and four tuples $(n,k)$ where $n$ is a leaf node in the subtree and $k$ is the predicted class at node $n$ in the optimal subtree.
\end{itemize}

We note that there is no guarantee of a unique optimal subtree, only a unique optimal objective value. The subset of samples correctly classified, the structure of the subtree and the number of leaves may differ between optimal solutions.

A basic outline is as follows. Given some integral solution\footnote{The MP constraints ensure that the tree structure is valid, but we do not check that the solution is feasible with respect to the full model. Our experiments suggest that it is advantageous to run the primal heuristic even if the solution is cut off by the Benders cuts.} $(b^*,p^*,w^*,\theta^*)$ corresponding to tree $\mathcal{DT}(b^*,p^*,w^*)$ that we wish to further optimise, we choose some set of internal nodes to be the roots of the subtrees which will be optimised. For each subtree root node we use D2S to find the optimal subtree, if the updated tree containing the optimised subtrees is better than the incumbent solution then insert the updated tree as the new incumbent. As subtree roots we choose nodes $n_r$ for which two conditions hold:

\begin{enumerate}
    \item The height of $n_r$ is at most $2$ and,
    \item Either the height of $n_r$'s parent node is greater than $2$ or $n_r$ is the root node.
\end{enumerate}

The first condition ensures that the subtree roots are not too far away from leaf nodes. The intention is that we optimise the tails of $\mathcal{DT}(b^*,p^*,w^*)$ without extending further than needed into nodes which are cut off in the solution. The second condition checks that the parent node isn't a more suitable candidate, the root node having no parent is treated as a special case. Denoting the subset of samples which pass through node $n$ as $\mathcal{I}_n(\mathcal{DT}(b^*,p^*,w^*))$, this procedure is formally described in Algorithm \ref{alg:Solution Polishing Primal Heuristic}.

\begin{algorithm}
\caption{Solution Polishing Primal Heuristic}
\label{alg:Solution Polishing Primal Heuristic}
\hspace*{\algorithmicindent} \textbf{Input: } Integral MP solution $(b^*, p^*, w^*, \theta^*)$, Leaf penalty $\lambda$ 
    \begin{algorithmic}[1]
        \State $\mathcal{R} \gets \{n \in \mathcal{B} \cup \mathcal{T} : (height(n) \leq 2) \land (n=1 \vee height(a(n)) > 2) \}$ \Comment{Collect subtree roots} 
        \For{$n \in \mathcal{R}$}
            \State $\bar{\mathcal{I}} \gets \mathcal{I}_n(\mathcal{DT}(b^*,p^*,w^*))$
            \State $(\hat{\mathcal{I}}^+, \hat{\mathcal{I}}^-, \mathcal{B}_{sub}, \mathcal{L}_{sub}) \gets D2S(\bar{\mathcal{I}}, \lambda)$
            \If{$|\hat{\mathcal{I}}^+| - \lambda |\mathcal{I}| |\mathcal{L}_{sub}| > \sum_{i \in \bar{\mathcal{I}}}\theta_i^* - \lambda |\mathcal{I}| \sum_{n \in \mathcal{D}^+(n)} p_n^*$} 
                \State $\text{Update } \mathcal{DT}(b^*,p^*,w^*) \text{ with optimal subtree solution}$
            \EndIf
        \EndFor 

    \Return $\mathcal{DT}(b^*,p^*,w^*)$
    \end{algorithmic}
\end{algorithm}

In practice the same subtrees are seen and re-optimised repeatedly. To avoid recomputation the optimal subtree solutions are cached, keyed on the path through the upper portion of the tree which leads into the subtree.

\subsection{Equivalent Point Valid Inequalities}\label{sect: Equivalent Point Valid Inequalities}
In this section we detail a class of inequalities which exploits the fact that samples with identical or similar feature values often have different classes. If it can be shown that such samples will always fall into the same leaf for a given tree then a bound can be established on the classification score of those samples, since we would never be able to simultaneously classify all of those samples correctly. This idea has been exploited to great effect in DP formulations to tighten bounds which prune the search space  \citep{huOptimalSparseDecision2019, linGeneralizedScalableOptimal2020} and to compress the dataset size \citep{demirovicBlossomAnytimeAlgorithm2023}. Modelling equivalent points is complicated in our MIP setting since it requires encoding the aforementioned logic as linear constraints, but it also allows us to model \textit{almost} equivalent points, those that are equivalent in all but a small subset of features.

To model this, we introduce the idea of equivalent point sets (EQP sets) which are sets of samples which will be routed into the same leaf given that some condition on the tree structure holds. We say that a subset $J \subseteq I$ is an EQP set with respect to a set of features $\bar{F} \subseteq F$ if the following conditions hold:

\begin{equation}\label{EQP Set Definition}
    \begin{gathered}
        x_f^{i_1} = x_f^{i_2},\quad \forall f \in \bar{F},\; \forall (i_1,i_2) \in J \times J \\ 
        \exists (i_1,i_2) \in J \times J :y^{i_1} \neq y^{i_2}
    \end{gathered}
\end{equation}

That is to say that the samples in $J$ agree for all feature values in $\bar{F}$ but do not all share the same class. For each EQP set $J$ there is an associated set of features $F^*$ for which the feature values do not agree, which is to say that $J$ is an EQP set with respect to $F \setminus F^*$. We refer to $F^*$ as a split set, owing to its ability to split the samples in $J$ when branched on. For example, consider three samples with feature values $x_1 = (1,0,1)$, $x_2 = (1,0,1)$, $x_3 = (1,0,0)$ and $y_1 = 0$, $y_2 = 1$, and $y_3 = 2$. $J = \{1,2,3\}$ is an EQP set with respect to $\bar{F} = \{0,1\}$ since the samples have different classes and agree in all but the last feature. If the samples are not routed through a node which branches on the last feature then at most one out of the three samples can be correctly classified.

For a given EQP set $J$ with associated split set $F^*$ we define two constraint sets which together define valid inequalities. The first set of constraints $\mathcal{H}(J,F^*)$ models the relationship between the branch variables $b$ and the routing of samples in $J$. The constraints in $\mathcal{H}$ detect when the tree structure guarantees that the samples in $J$ will be routed into the same leaf node. The second set of constraints $\mathcal{G}(J)$ bounds the classification score of the samples in the EQP set. A new set of decision variables $\beta$ is introduced to mediate between the two constraint sets. Its particular structure depends on the how $\mathcal{H}$ is defined but for the variants of $\mathcal{H}$ and $\mathcal{G}$ which we introduce the following must hold:

\begin{enumerate}
    \item It requires some special element $\beta_J^\mathcal{G}$.
    \item $\beta_J^\mathcal{G}$ is forced to zero whenever $b$ is such that we can guarantee that the samples in $J$ flow into the same leaf node and is otherwise unbounded.
\end{enumerate}

Given these constraint sets, we add the following constraints to the MP (\ref{model:BendOCT}) for each EQP set:

\begin{align}
    &(b, \beta) \in \mathcal{H}(J, F^*) \label{const: EQP Linking Constraints}\\
    &(\theta, \beta) \in \mathcal{G}(J) \label{const: EQP Bounding Constraints}
\end{align}

In the variants of $\mathcal{H}$ and $\mathcal{G}$ which we introduce, $\mathcal{H}$ will force $\beta_J^\mathcal{G}$ to zero whenever it can be guaranteed that the samples in $J$ flow into the same leaf, which will enforce a bound on the classification of samples in $J$ via the constraints in $\mathcal{G}$. When such a guarantee cannot be made, $\mathcal{H}$ will leave $\beta_J^\mathcal{G}$ unbounded, allowing $\mathcal{G}$ to relax the bound on sample classification.

In general there are $\mathcal{O}(|\mathcal{I}|^2)$ EQP sets which could be added to the model. Modelling all possible EQP sets is counterproductive since the number of variables and constraints added slows down solving the LP relaxation without proportionately tightening the relaxation. The choice of which EQP sets to add is discussed in Section \ref{Generating Equivalent Point Sets}.

\subsubsection{Basic EQP Cuts}

We first consider basic versions of $\mathcal{H}$ and $\mathcal{G}$. The simplest method to link the tree structure to the sample paths is to note that if none of the features in $F^*$ which split $J$ are branched on, then the samples in $J$ will necessarily flow into the same leaf node. We encode this logic as follows:

\begin{align}
    \mathcal{H}^{Basic}(J, F^*) = 
    \left\{
        \begin{aligned}
            (b,\beta) : &\beta_J^\mathcal{G} \leq \sum_{\substack{n \in \mathcal{B} \\ f \in F^*}} b_{nf} \\
            &\beta_J^\mathcal{G} \in [0,1]
        \end{aligned}
     \right\} 
\end{align}

If the tree branches on any of the features in $F^*$ then it can no longer be guaranteed that the samples in $J$ fall into the same leaf and $\beta_J^\mathcal{G}$ is unbounded. Otherwise, $\beta_J^\mathcal{G}$ is forced to zero. The simplest method to bound the classification of samples in $J$ is to define $m(J) = \max_{k \in \mathcal{K}}\left| \{ i \in J : y^i = k \} \right|$ and enforce an upper bound of $m(J)$ whenever $\beta_J^\mathcal{G} = 0$:

\begin{equation}
    \mathcal{G}^{Basic}(J) =  \left\{ (\theta, \beta) : \sum_{i \in J} \theta_i \leq m(J) + \left[ |J| - m(J) \right] \beta_J^\mathcal{G} \right\}
\end{equation}

We stress that $\mathcal{H}^{Basic}$ never forces $\beta_J^\mathcal{G}$ to be equal to one. Instead, $\mathcal{G}^{Basic}$ implements a bound on some subset of the objective. As such the objective exerts an upwards pressure on $\beta_J^\mathcal{G}$ which can force it to one as necessary. 

\subsubsection{Chain EQP Cuts}

It is possible for the tree to branch on some feature $f \in F^*$ at node $n$ without the samples in $J$ actually passing through node $n$. In this case $\mathcal{H}^{Basic}$ is overly conservative since the bound in $\mathcal{G}$ will be relaxed even though the samples actually flowed into the same leaf. Remedying this requires modelling the paths of the EQP sets through the tree to detect when the samples are actually split. We say that EQP set $J$ is split by $F^*$ at internal node $n \in \mathcal{B}$ if the samples in $J$ follow a path from the root node to node $n$ and node $n$ branches on a feature in $F^*$ (i.e. $\sum_{f \in F^*} b_{nf} = 1$). We introduce binary decision variables $\beta_J^n$ which are equal to one only if $J$ is split by $F^*$ at node $n \in \mathcal{B}$. Then, the sample paths are modelled as follows:

\begin{equation}
    \mathcal{H}^{Chain}(J, F^*) = 
    \left\{
        \begin{aligned}
            (b,\beta) : &\beta_J^n \leq \frac{1}{|\mathcal{A}(n)| + 1} \left[ \begin{aligned}
                 \sum_{f \in F^*}b_{nf} + &\sum_{\substack{n_a \in \mathcal{A}^L(n) \\ f \in F\setminus F^* : x_f^j = 0}} b_{n_af} \\ +&\sum_{\substack{n_a \in \mathcal{A}^R(n) \\ f \in F\setminus F^*: x_f^j = 1}} b_{n_af} 
            \end{aligned} \right], \quad \forall n \in \mathcal{B} \\
            &\beta_J^\mathcal{G} \leq \sum_{n \in \mathcal{B}} \beta_J^n \\
            &\beta_J^n \in \{ 0, 1\},\quad \forall n \in \mathcal{B} \\
            &\beta_J^\mathcal{G} \in [0,1]
        \end{aligned}
     \right\} 
\end{equation}

The first constraint encodes that $\beta_J^n$ is forced to zero unless all of the following hold:
\begin{enumerate}
    \item Node $n$ branches on some feature $f \in F^*$.
    \item At all ancestors of $n$, the tree branches on a feature $f \in F \setminus F^*$ such that all samples in $J$ remain on the path to node $n$. i.e. they are not split or branch into a subtree which does not contain $n$.
\end{enumerate}

When all of these conditions hold the constraint reads $\beta_J^n \leq 1$. If any of the conditions do not hold then $\beta_J^n < 1$, which forces $\beta_J^n$ to zero since it is a binary variable. The second constraint ensures that $\beta_J^\mathcal{G}$ is forced to zero if $J$ is not split at any of the nodes in the tree. The first constraint can be disaggregated as follows:

\begin{equation}
    \mathcal{H}_{DA}^{Chain}(J, F^*) = 
    \left\{
        \begin{aligned}
            (b,\beta) : &\beta_J^n \leq \sum_{f \in F^*}b_{nf},  \quad \forall n \in \mathcal{B} \\ 
            &\beta_J^n \leq \sum_{f \in F\setminus F^* : x_f^j = 0} b_{n_af}, \quad  \forall n \in \mathcal{B}, \forall n_a \in \mathcal{A}^L(n) \\
            &\beta_J^n \leq \sum_{f \in F\setminus F^* : x_f^j = 1} b_{n_af}, \quad  \forall n \in \mathcal{B}, \forall n_a \in \mathcal{A}^R(n) \\
            &\beta_J^\mathcal{G} \le \sum_{n \in \mathcal{B}} \beta_J^n \\
            &\beta_J^n \in [0,1], \quad \forall n \in \mathcal{B} \\
            &\beta_J^\mathcal{G} \in [0,1]
        \end{aligned}
     \right\} 
\end{equation}

Disaggregation requires many more constraints but should have a stronger relaxation (and indeed does experimentally) and also allows for the integrality of $\beta$ to be relaxed. 

\subsubsection{Recursive EQP Cuts}
An alternative method of modelling sample paths is to define them recursively. In this case we take $\beta_J^n \in [0,1]$ to be strictly equal to zero unless there exists a path from node $n$ to some descendant $n^*$ which splits on $F^*$ (we also consider it a path if $n = n^*$). That is to say that 
$\beta_J^n$ is unbounded if any of the following hold:
\begin{itemize}
    \item Node $n$ branches on some feature $f \in F^*$.
    \item $\beta_J^{l(n)} = 1$ and all samples in $J$ branch into node $l(n)$.
    \item $\beta_J^{r(n)} = 1$ and all samples in $J$ branch into node $r(n)$.
\end{itemize}

We introduce auxiliary variables $\beta_J^{n,l} \in \{ 0,1 \}$ at each branch node $n$ which are equal to one if the samples in $J$ branch into the left child of $n$ and $\beta_J^{l(n)} = 1$. $\beta_J^{n,r}$ is defined analogously. Then the sample paths can be modelled recursively as follows:

\begin{equation}
    \mathcal{H}^{Recursive}(J, F^*) = 
    \left\{
        \begin{aligned}
            (b,\beta) : &\beta_J^n \leq \sum_{f \in F^*}b_{nf} + \beta_J^{n,l} + \beta_J^{n,r}, \quad \forall n \in \mathcal{B} \\
            &\beta_J^{n,l} \leq \frac{1}{2} \left[ \sum_{f \in F \setminus F^* : x_f^i=0} b_{nf} + \beta_J^{l(n)} \right],\quad \forall n \in \mathcal{B} \\
            &\beta_J^{n,r} \leq \frac{1}{2} \left[ \sum_{f \in F \setminus F^* : x_f^i=1} b_{nf} + \beta_J^{r(n)} \right],\quad \forall n \in \mathcal{B} \\
            &\beta_J^\mathcal{G} \leq \beta_J^1 \\
            &\beta_J^n \in [0,1], \quad \forall n \in \mathcal{B} \\
            &\beta_J^{n,l} \in \{ 0,1\}, \quad \forall n \in \mathcal{B} \\
            &\beta_J^{n,r} \in \{ 0,1\}, \quad \forall n \in \mathcal{B} \\
            &\beta_J^\mathcal{G} \geq 0
        \end{aligned}
     \right\} 
\end{equation}

Some edge cases have been omitted for clarity. In practice the constraints must be modified to account for the fact that the children of some internal nodes are terminal nodes and do not have associated $\beta$ variables. Just as for $\mathcal{H}^{Chain}$ the constraints defining $\beta_J^{n,l}$ and $\beta_J^{n,r}$ can be disaggregated. This requires two extra constraints per internal node but allows for integrality to be relaxed on all $\beta$ variables.

\begin{equation}
    \mathcal{H}_{DA}^{Recursive}(J, F^*) = 
    \left\{
        \begin{aligned}
            (b,\beta) : &\beta_J^n \leq \sum_{f \in F^*}b_{nf} + \beta_J^{n,l} + \beta_J^{n,r}, \quad \forall n \in \mathcal{B} \\
            &\beta_J^{n,l} \leq  \sum_{f \in F \setminus F^* : x_f^i=0} b_{nf}, \quad \forall n \in \mathcal{B} \\
            &\beta_J^{n,r} \leq  \sum_{f \in F \setminus F^* : x_f^i=1} b_{nf}, \quad \forall n \in \mathcal{B} \\
            &\beta_J^{n,l} \leq \beta_J^{l(n)}, \quad \forall n \in \mathcal{B} \\
            &\beta_J^{n,r} \leq \beta_J^{r(n)}, \quad \forall n \in \mathcal{B} \\
            &\beta_J^\mathcal{G} \leq \beta_J^1 \\
            &\beta_J^n \in [0,1], \quad \forall n \in \mathcal{B} \\
            &\beta_J^{n,l} \in [0,1], \quad \forall n \in \mathcal{B} \\
            &\beta_J^{n,r} \in [0,1], \quad \forall n \in \mathcal{B} \\
            &\beta_J^\mathcal{G} \geq 0
        \end{aligned}
     \right\} 
\end{equation}

\subsubsection{Group Selection Constraints}
Given that the samples in an EQP set $J$ path into the same leaf as indicated by $\beta_J^\mathcal{G}$ being forced to zero, $\mathcal{H}^{Basic}$ bounds the classification score of the samples in $J$ based on the majority class. An alternative approach is to partition $J$ into groups based on the sample classes and enforce that only samples from one of these groups can be correctly classified at a time. We index the groups by $\mathcal{K}_J^* = \{k \in \mathcal{K} : \left( \exists i \in J : y_i = k \right) \}$ and denote the samples in group $k \in \mathcal{K}_J^*$ by $J_k$. For each group we define decision variable $G_k \in [0,1]$ which is equal to one if the samples in group $k \in \mathcal{K}_J^*$ are free to be correctly classified and zero otherwise. The group selection constraints are then:

\begin{equation}
    \mathcal{G}^{GS}(J) =  
    \left\{  
        \begin{aligned}
            (\theta, \beta) : & \sum_{k \in \mathcal{K}_J^*}G_k \leq 1 + \left( |\mathcal{K}_J^*| - 1 \right) \beta_J^\mathcal{G} \\
            &\sum_{i \in J_k} \theta_i = |J_k|G_k, \quad \forall k \in \mathcal{K}_J^* \\
            &G_k \in [0,1], \quad \forall k \in \mathcal{K}_J^*
        \end{aligned}
    \right\}
\end{equation}

When $\beta_J^\mathcal{G} = 0$ the first constraint asserts that at most one of the group selection variables $G_k$ can be equal to one, and the second constraint enforces that only the group with $G_k = 1$ can have correctly classified samples. If $\beta_J^\mathcal{G}$ is not forced to zero by $\mathcal{H}$ then the first constraint is inactive and the $G_k$'s vary freely to match $\theta$ as per the second constraint.

\subsubsection{Generating Equivalent Point Sets}\label{Generating Equivalent Point Sets}
Thus far we have assumed that we are given the EQP sets to be modelled, but in general it is not obvious which EQP sets to model. As defined by conditions (\ref{EQP Set Definition}) there will be $\mathcal{O}(|\mathcal{I}|^2)$ EQP sets, including all of them may produce a tight relaxation but will make the relaxation far too large to solve efficiently. In this work we include EQP sets for which the associated split set $F^*$ is small. Experimentally we find $|F^*| \leq 2$ is a reasonable bound in most cases. The intuition is that larger split sets produce more opportunities for the bound to be relaxed, reducing the marginal benefit of modelling the associated EQP set. This is particularly true for variants of $\mathcal{H}$ which model the sample paths since they require many extra variables and constraints.

Another choice remains of when to aggregate EQP sets. We choose to aggregate any EQP sets which share a split set. For example EQP sets $J_1 = \{1, 2 \}$ and $J_2 = \{2,3 \}$ with shared split set $F^*$ would be combined into one EQP set $J = \{1,2,3\}$. In many cases the constraints in $\mathcal{G}^{Basic}$ can be tightened under the intuition that after an EQP set is split, the resulting subsets may themselves be EQP sets. This can also be trivially deduced from the constraints themselves and experiments indicated that modern MIP solvers could detect the tightening during presolve. As such we make no further modifications to the EQP sets.

\subsection{Path Bound Cutting Planes}\label{sect:Path Bound Cutting Planes}

\subsubsection{Basic Cutting Planes}

\begin{figure}
    \centering
    \includegraphics[width=0.5\linewidth]{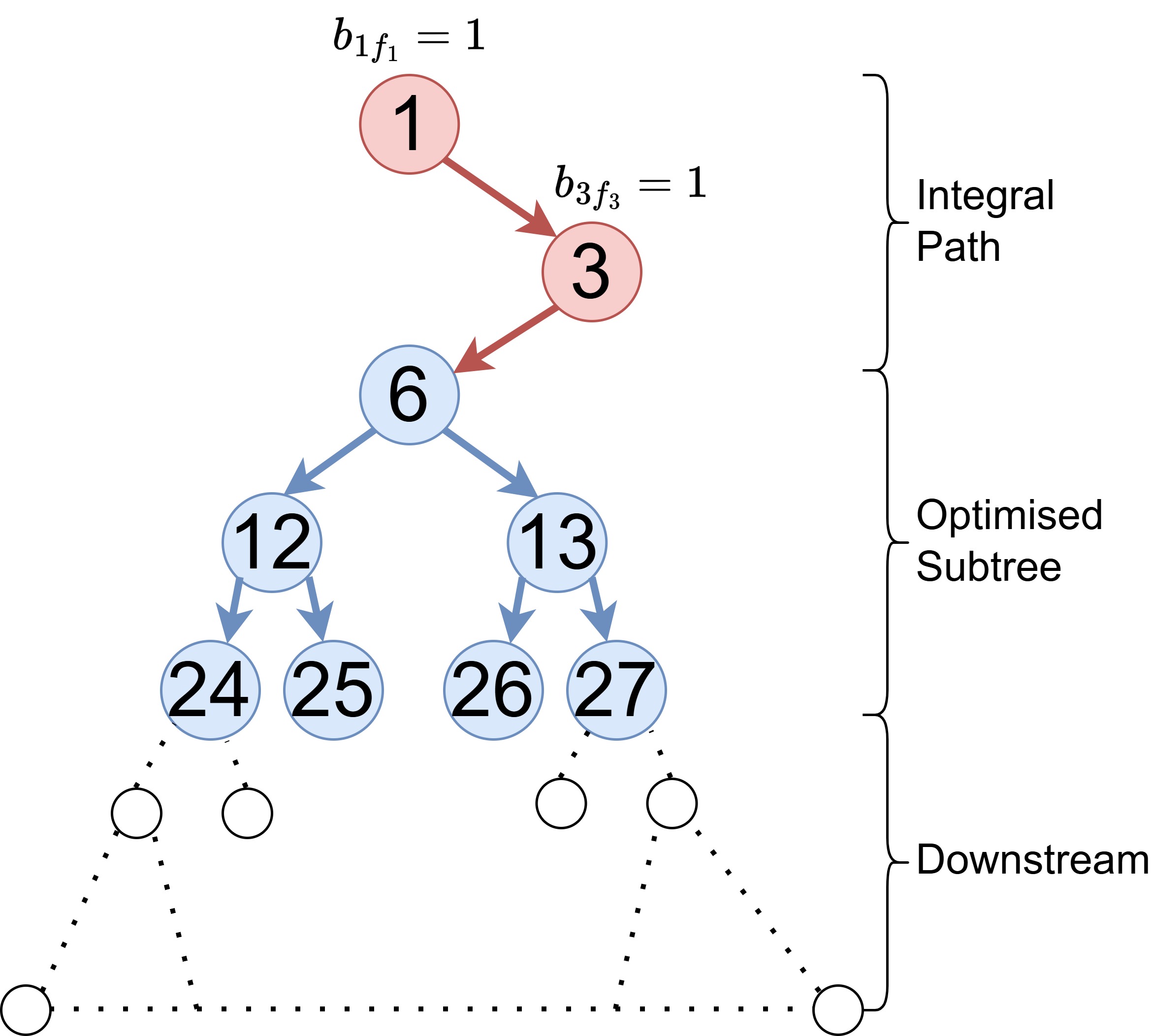}
    \caption{A subset of a tree with an integral path in the relaxation solution. It is partitioned into the integral path, the depth two subtree to be optimised, and the nodes downstream of the optimised subtree. Decision variables associated with nodes not in the integral path may be fractional.}
    \label{fig:PBCP Diagram}
\end{figure}

The D2S subroutine introduced in Section \ref{sect:Solution Polishing Primal Heuristic} and detailed in \ref{Appendix:Optimal Depth 2 Tree Subroutine} efficiently finds optimal depth two subtrees. In the solution polishing primal heuristic we fix the upper layers of trees and use the D2S subroutine to optimise the tails, but the optimality of the subtree provides further information that is being discarded. The returned subtree also implies an upper bound on the classification of samples which are routed into the subtree whenever the path into that subtree is present in the tree and the leaf nodes in the subtree are limited to a depth of at most two. That is to say that under the conditions which the optimal subtree was generated we cannot have a better classifier of the samples which were routed into the subtree. We introduce path bound cutting planes to exploit this. To model this we introduce the following notation:

\begin{itemize}
    \item Denote a path through a tree of length $m$ as $P = (r_1,...,r_m)$. The elements $r$ of the path are 3-tuples $(n,f,dir)$ where $n$ is a node $n \in \mathcal{B}$, $f$ is the branch feature at node $n$, and $dir \in \{0,1 \}$ indicates the direction (0 for left, 1 for right) the path follows from node $n$. A path $P$ terminates at some node not in $P$ itself, E.g. The final element of the integral path in Figure \ref{fig:PBCP Diagram} is $(3, f_3, 0)$ and terminates at node $6$.
    \item Denote $\mathcal{I}(P) = \{i \in I : x_f^i = dir, \; \forall (.,f,dir) \in P \}$ to be the subset of $\mathcal{I}$ which follows path $P$ through the tree.
    \item Define an integral path $P$ through tree $\mathcal{DT}(\bar{b},\bar{p},\bar{w})$ to be one such that $b$ takes integral values at all nodes in the path, i.e. $\bar{b}_{nf} = 1 \; \forall (n,f,.) \in P$. There is a subset of samples which follows an integral path without being "split" by a fractional branching rule, even though the relaxation solution which defines the tree may itself be fractional.
    \item We say that a path $P$ terminates at a node $n_{sub}$ which is the root of a subtree with nodes $\mathcal{D}(P) = \mathcal{D}^+(n_{sub})$. As shown in Figure \ref{fig:PBCP Diagram} we partition $\mathcal{D}(P)$ into the nodes in the depth two subtree which will be optimised, $\mathcal{D}^{sub}(P) = \{n_d \in  \mathcal{D}(P) : dist(n_{sub},n_d) \leq 2\}$, and the nodes which are downstream of the optimised subtree, $\mathcal{D}^{down}(P) = \{n_d \in  \mathcal{D}(P) : dist(n_{sub},n_d) > 2\}$.
\end{itemize}

The cuts are implemented as follows. At non-integral nodes in the branch and bound tree, retrieve possibly fractional relaxation solution $(\bar{b}, \bar{p}, \bar{w}, \bar{\theta})$. Search for integral paths in the tree $\mathcal{DT}(\bar{b}, \bar{p}, \bar{w})$ which terminate at a node with a height of at least two from the terminal nodes $\mathcal{T}$ (not the height in the tree $\mathcal{DT}(\bar{b}, \bar{p}, \bar{w})$). This includes paths which may be subsets of other integral paths. Whenever such an integral path P is found, calculate $(\hat{\mathcal{I}}^+, \hat{\mathcal{I}}^-, \mathcal{B}_{sub}, \mathcal{L}_{sub}) = D2S(\mathcal{I}(P))$ which provides information\footnote{See Section \ref{sect:Solution Polishing Primal Heuristic} for a full description of the output of D2S} about the optimal depth two subtree for samples following path $P$. If the current relaxation solution in the subtree is better than the optimal subtree solution then we add a cut to the model. When the cut is active it should bound the total classification score of samples in $\mathcal{I}(P)$ by $|\hat{\mathcal{I}}^+|$ while relaxing the bound under the following conditions:

\begin{itemize}
    \item A different branch decision is made along path $P$.
    \item A prediction is made somewhere along path. $P$
    \item One of the nodes downstream of the optimised subtree is a leaf node.
\end{itemize}

To encode these conditions we set

\begin{equation}
    {relax}^P = \sum_{(n,f^*,.) \in P} \left(p_n + \sum_{f \in F : f \neq f^*}b_{nf}\right)  + \sum_{n \in \mathcal{D}^{down}(P)} p_n 
\end{equation}

after which the cut follows:

\begin{equation}\label{const:Basic Path Bound Cutting Planes Cut}
    \sum_{i \in \mathcal{I}(P)} \theta_i \leq |\hat{\mathcal{I}}^+| + [|\mathcal{I}(P)| - |\hat{\mathcal{I}}^+|] * {relax}^P 
\end{equation}

\begin{algorithm}[h]
\caption{Path Bound Cutting Planes}
\label{alg: Path Bound Cutting Planes}
\hspace*{\algorithmicindent} \textbf{Input: } Fractional MP solution $(\bar{b}, \bar{w}, \bar{\theta})$, Leaf penalty $\lambda$ 
    \begin{algorithmic}[1]
        \State $P^{int} \gets \emptyset$ \Comment{Initialise set of integral paths}
        \State $S \gets \{(1, [\;] )\}$ \Comment{Initialise stack with root node and empty path}
        \While{$S \neq \emptyset$}
            \State Pop $(n, P^{partial})$ from $S$
            \If{$height(n) > 2$}
                \For{$f \in \mathcal{F}$}
                    \If{$\bar{b}_{nf} = 1$}
                        \State Create $P_{left}^{partial}$ by appending $(n, f, 0)$ to $P^{partial}$
                        \State Create $P_{right}^{partial}$ by appending $(n, f, 1)$ to $P^{partial}$
                        \State Add $(l(n), P_{left}^{partial})$ to $S$
                        \State Add $(r(n), P_{right}^{partial})$ to $S$
                        \State Add $P_{left}^{partial}$ and $P_{right}^{partial}$ to $P^{int}$
                        \State \algorithmicbreak{}
                    \EndIf
                \EndFor
            \EndIf
        \EndWhile
        \For{$P \in P^{int}$}
            \State $(\hat{\mathcal{I}}^+, \hat{\mathcal{I}}^-, \mathcal{B}_{sub}, \mathcal{L}_{sub}) \gets D2S(\mathcal{I}(P), \lambda)$
            \If{$\sum_{i \in \bar{\mathcal{I}}(P)}\bar{\theta}_i - \lambda |\mathcal{I}| \sum_{n \in \mathcal{D}(P)} \bar{p}_n >|\hat{\mathcal{I}}^+| - \lambda |\mathcal{I}| |\mathcal{L}_{sub}|$}
                \State Add constraint (\ref{const:Basic Path Bound Cutting Planes Cut}) to MP
            \EndIf
        \EndFor
    \end{algorithmic}
\end{algorithm}

We describe this procedure formally in Algorithm \ref{alg: Path Bound Cutting Planes}. In practice the same integral paths are seen repeatedly in relaxation solutions. To avoid recomputation optimal subtree solutions are cached. The cache design must take into account that different paths may filter out the same subset of the samples, i.e. $\mathcal{I}(P_1) = \mathcal{I}(P_2)$ for distinct paths $P_1$ and $P_2$. For example a path which branches right on feature $f_1$ then left on feature $f_2$ is effectively identical to a path which branches left on feature $f_1$ then right on feature $f_2$. The cache is keyed on the set of features and directions branched on in the path to capture this symmetry.

\subsubsection{Improved Cutting Planes}

Experimentally the cuts are effective as described, but we note two avenues for improvement:

\begin{enumerate}
    \item These cuts often produce a situation in which the solver has a tight bound but is slow to converge to an optimal solution, even though the D2S subroutine has likely found subtrees which are present in optimal solutions. Passing the optimised subtrees back to the solver as suggested solutions where possible should accelerate convergence.
    \item The D2S subroutine returns much more information than is made use of, in particular the structure of the optimal subtree and the classification score of specific samples. Exploiting this additional information should prevent the solver from branching through symmetric optimal subtrees and yield stronger cuts.
\end{enumerate}

One solution to the first item is to apply the cuts in tandem with the solution polishing primal heuristic described in Section \ref{sect:Solution Polishing Primal Heuristic}. The cache of optimal subtrees which have been computed can be shared to minimise recomputation and provides an avenue for optimised subtrees to be passed back to the solver. On the second item we propose two extensions of the cutting planes which exploit the extra information. The first alternative, which we refer to as \textit{bounding negative samples}, is a modification of cut (\ref{const:Basic Path Bound Cutting Planes Cut}) which forces $\theta_i = 0$ for all samples which are misclassified in the optimal subtree whenever $relax^P=0$, instead of enforcing a bound over all samples which follow the path. The negative samples are bounded by the following cut:

\begin{equation}\label{const:Bound Negative Samples Cut}
    \sum_{i \in \hat{\mathcal{I}}^-} \theta_i \leq |\hat{\mathcal{I}}^-| * {relax}^P
\end{equation}

The second modification is to \textit{bound the tree structure} by restricting the decision variables corresponding to the structure of the subtree to agree with that of the optimal subtree whenever ${relax}^P = 0$. This is implemented by forcing all of the branch and prediction variables not in $\mathcal{B}_{sub}$ or $\mathcal{L}_{sub}$ to zero:

\begin{equation}
    \sum_{(n,f^*) \in \mathcal{B}_{sub}} \left(p_n + \sum_{f \in F : f \neq f^*}b_{nf}\right) + \sum_{(n,k^*) \in \mathcal{L}_{sub}} \left( \sum_{f \in F} b_{nf} + \sum_{k \in K : k \neq k^*}w_k^n \right) \leq (|\mathcal{B}_{sub}| + |\mathcal{L}_{sub}|) * {relax}^P
\end{equation}

\section{Results}\label{sect: Computational Experiments}

This section details the experiments run to assess the effectiveness of the proposed acceleration techniques. Since we are interested in accelerating convergence we compare how many instances can be solved to optimality and the optimality gap in those that can not. We do not compare out-of-sample performance.

\subsection{Binary Data Encoding}\label{Binary Data Encoding}

\begin{figure}[h]
    \centering
    \includegraphics[width=0.8\linewidth]{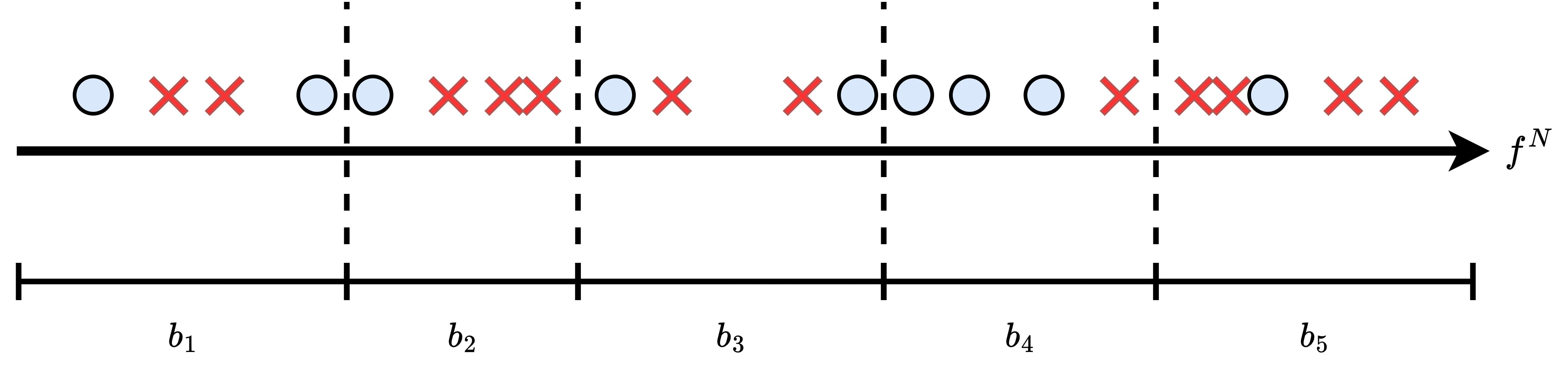}
    \caption{Example quantile bucket encoding on continuous feature $f^N$ with two classes. $f_N$ will have five associated binary features where the $m^{th}$ binary feature is equal to one whenever the feature value falls into the $m^{th}$ bucket and zero otherwise.}
    \label{fig:Bucket Encoding}
\end{figure}

\begin{figure}[h]
    \centering
    \includegraphics[width=0.8\linewidth]{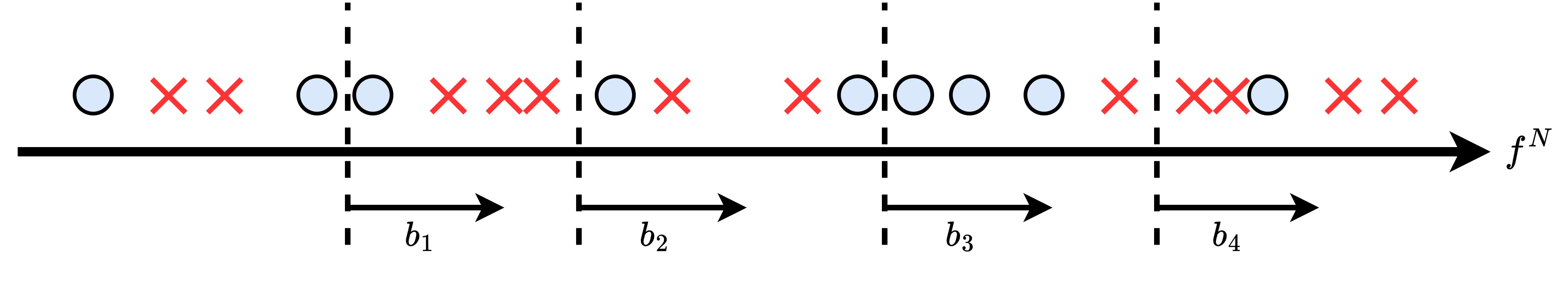}
    \caption{Example quantile threshold encoding on continuous feature $f^N$ with two classes. $f_N$ will have four associated binary features where the $m^{th}$ binary feature is equal to one whenever the feature value is greater than or equal to the $m^{th}$ threshold and zero otherwise.}
    \label{fig:Threshold Encoding}
\end{figure}

The BendOCT model requires datasets with binary features. This restriction is common, though not ubiquitous, in methods for constructing optimal classification trees. In the likely case that the dataset of interest has non-binary features it is necessary to encode those features as binary values. For categorical features we take the standard approach of one-hot encoding. Ordinal variables with many levels as treated as continuous variables. For datasets with continuous features we run experiments using two common methods of binarisations by bucketing or thresholding continuous feature values as shown in Figures \ref{fig:Bucket Encoding} and \ref{fig:Threshold Encoding}. In this work the position of the buckets and the thresholds corresponds to the 5-quantiles of the data. As such we refer to the binarisation schemes as QB-5 and QT-5.

\subsection{Experimental Setup}
All experiments were run on an AMD EPYC 9684X processor with 4GB of RAM. We used Python 3.12.7 to code all experiments, Gurobi 12.0.0 \citep{gurobioptimizationllcGurobiOptimizerReference2024} was used as a MIP solver and the D2S subroutine was implemented using Numpy 2.2.4. All experiments were limited to a single thread for benchmarking purposes with a time limit of one hour and an optimality gap of $0\%$. All other Gurobi parameters were left at the default values. A total of 33 datasets from the UCI machine learning repository \citep{kellyUCIMachineLearningn.d.} are considered. Of these 12 have purely categorical features, 12 purely continuous features, and the remaining 9 a mix of categorical and continuous features. Categorical features are one-hot encoded, continuous features are discretised by either quantile bucketing (QB-5) or quantile thresholding (QT-5). The datasets, range between $47$ and $7200$ samples, and between $15$ and $357$ binary features. Table \ref{tab: Dataset Statistics} in the Appendix list all datasets used.

Experiments are run over 1620 unique instances wherein each instance constitutes a unique combination of the dataset, maximum tree depth $D \in \{3,4\}$, regularisation parameter \newline $\lambda \in \{0.08, 0.06, 0.04, 0.02, 0.01,$ $0.008, 0.006, 0.004, 0.002, 0.001, 0.0008, 0.0006, 0.0004, 0.0002, 0.0001\}$, and binarisation scheme for datasets with continuous features. We warm start the solver with a CART heuristic trained on binarised features, choosing splitting rules by minimising Gini impurity and pruning the tree by cost-complexity pruning based on accuracy with leaf penalty $\lambda$. The initial solution is further improved by solution polishing as described in Section \ref{sect:Solution Polishing Primal Heuristic}. We use the BendOCT model as described in Section \ref{sect:BendOCT} as a baseline for comparison. Code to run all experiments and per instance results are available at https://github.com/mitchellkeegan/accel-OCT.


\subsection{Empirical Results}\label{sect: Main Results}
We report results for the best performing variants of each of the proposed methods. Extended results are available in \ref{Appendix: Extended Results} which assess the effectiveness of variants of the EQP inequalities and path bound cutting planes to justify the cut configurations which we use. We find that EQP inequalities are most effective with the disaggregated recursive variant of $\mathcal{H}$, group selection constraints, and a bound on the size of split sets of $|F^*| \leq 2$. Path bound cutting planes were found most effective when bounding negative samples and bounding the tree structure.

\begin{figure}[htb!]
    \centering
    \begin{subfigure}[t]{0.45\textwidth}
        \centering
        \includegraphics[width=\textwidth]{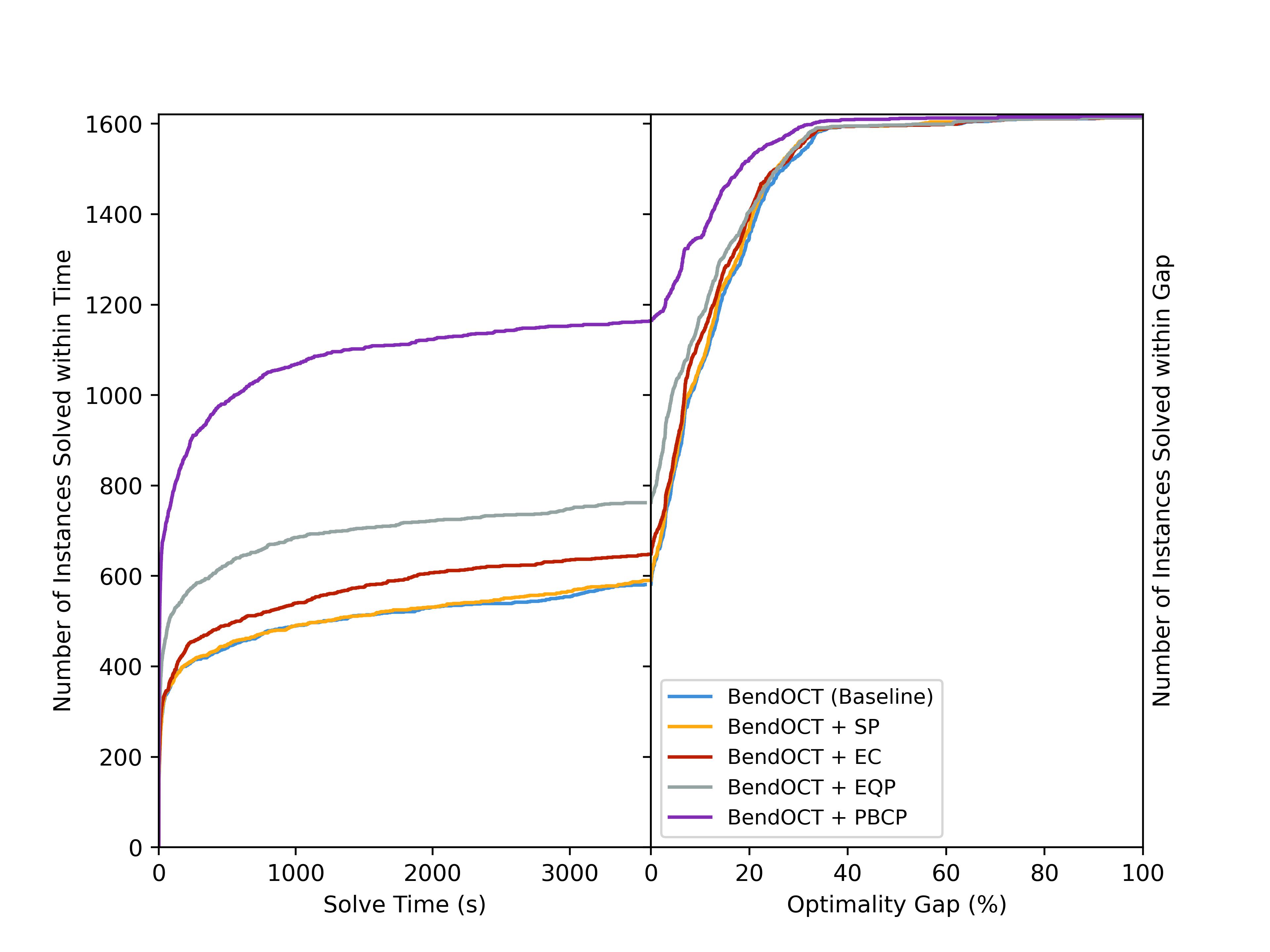}
        \caption{All Methods in Isolation}
        \label{fig: All Method Results}
    \end{subfigure}
    \begin{subfigure}[t]{0.45\textwidth}
        \centering
        \includegraphics[width=\textwidth]{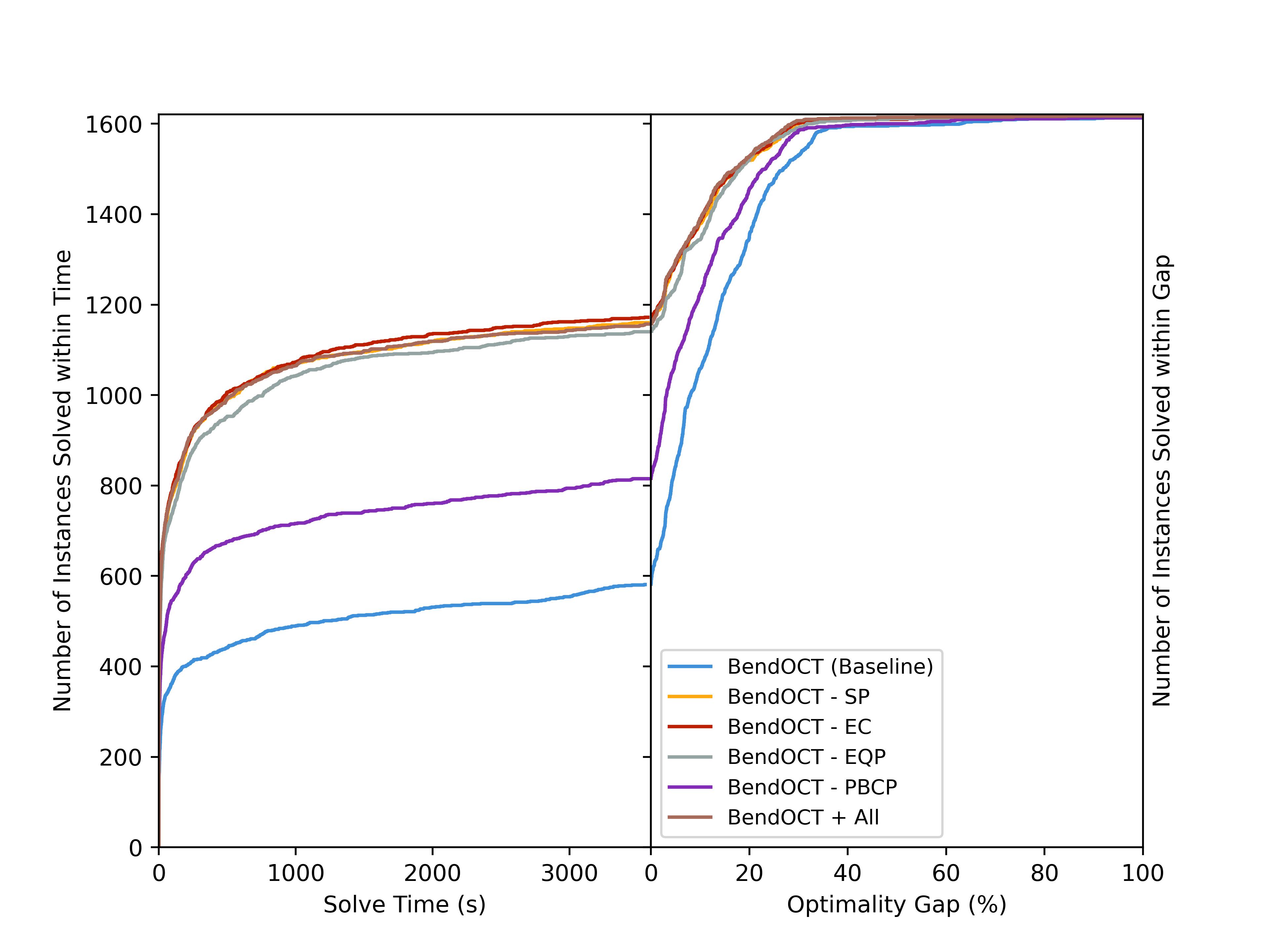}
        \caption{Ablation Test}
        \label{fig: Ablation Test Results}
    \end{subfigure}
    \caption{The left panels indicate how many instances each model solved to optimality within the given time. The right panels indicate for how many instances each model found a solution within the given optimality gap. We compare strengthened Benders cuts (EC), solution polishing (SP), path bound cutting planes (PBCP), and equivalent point bounds (EQP).}
    \label{fig:Main Results}
\end{figure}

We compare the performance of all methods implemented in isolation in Figure \ref{fig: All Method Results} and perform an ablation test in Figure \ref{fig: Ablation Test Results} to understand which techniques contribute the most towards performance. We find that in isolation the path bound cutting planes are by far the most effective, followed by EQP inequalities and strengthened Benders cuts. Solution polishing has a marginal effect, slightly improving the optimality gap in instances which cannot be solved to optimality. In the ablation test removing the path bound cutting planes has by far the biggest impact. Removing the EQP inequalities has a marginal negative impact and removing the solution polishing primal heuristic has little noticeable effect. Interestingly, the strengthened Benders cuts are mildly detrimental when combined with other acceleration techniques. Clearly the path bound cutting planes are the biggest driver of performance. This is despite requiring a call to the D2S subroutine during optimisation indicating that the information encoded in the cuts is extremely valuable to the solver. The vanilla BendOCT model is able to solve 582/1620 instances to optimality within the one hour time limit. Our best performing variant without the strengthened Benders cuts is able to solve as many instances within 12 seconds, a 300x speedup, and a total of 1173 instances within the time limit, over twice as many.

We summarise per dataset results in Table \ref{tab:Main Results Table}, aggregating over all available instances for each dataset. We compare BendOCT to accelerated BendOCT without strengthed Benders cuts. We also list results for a partially accelerated BendOCT with the strengthened cuts but without path bound cutting planes. This gives a more realistic idea of how the results may generalise to global objectives and constraints, since path bound cutting planes assume subtree independence, while other techniques are less restrictive. For example a hard constraint on the number of leaves would render the path bound cutting planes invalid as implemented currently, since we could no longer assert subtree optimality without knowing how many leaves have been used in the remainder of the tree. We find that both accelerated models dominate BendOCT, solving significantly more instances with shorter solve times and smaller optimality gaps in unsolved instances. The most extreme examples are the segmentation, spambase, and biodeg datasets for which BendOCT cannot solve any of the 60 instances, while the accelerated model solved 24, 20, and 29 instances respectively. The improvements are more modest for the partially accelerated model. This reinforces that path bound cutting planes are a large driver of performance, although we note that 10 out of 33 datasets have no EQP sets for the size of split sets we model while another four have relatively few which undermines the performance of the partially accelerated model.

\begin{table}[htb!]
    \centering
    \renewcommand*{\arraystretch}{0.8}
    \begin{tabular}{cccccccccc}
\hline
\multirow{2}{*}{Dataset} & \multicolumn{3}{c}{BendOCT} & \multicolumn{3}{c}{Accelerated BendOCT} & \multicolumn{3}{c}{Accelerated (No PBCP)}\\
  & Solved & Time & Gap & Solved & Time & Gap & Solved & Time & Gap\\
\hline
soybean-small & $30$ & $0.2$ & - & $\mathbf{30}$ & $\mathbf{0.1}$ & \textbf{-} & $30$ & $0.2$ & -\\
monk3 & $25$ & $258.5$ & $0.30$ & $\mathbf{30}$ & $\mathbf{6.7}$ & \textbf{-} & $30$ & $59.5$ & -\\
monk1 & $30$ & $6.4$ & - & $\mathbf{30}$ & $\mathbf{1.1}$ & \textbf{-} & $30$ & $1.3$ & -\\
hayes-roth & $20$ & $353.6$ & $2.72$ & $\mathbf{30}$ & $\mathbf{32.7}$ & \textbf{-} & $30$ & $80.7$ & -\\
monk2 & $15$ & $1642.6$ & $17.89$ & $\mathbf{30}$ & $\mathbf{146.9}$ & \textbf{-} & $20$ & $212.3$ & $7.80$\\
house-votes-84 & $26$ & $257.3$ & $0.27$ & $\mathbf{30}$ & $\mathbf{6.4}$ & \textbf{-} & $30$ & $16.5$ & -\\
spect & $18$ & $956.7$ & $3.63$ & $\mathbf{30}$ & $\mathbf{78.1}$ & \textbf{-} & $23$ & $230.2$ & $1.56$\\
breast-cancer & $7$ & $468.7$ & $16.26$ & $\mathbf{19}$ & $\mathbf{79.5}$ & $\mathbf{9.57}$ & $17$ & $1582.3$ & $9.10$\\
balance-scale & $\mathbf{18}$ & $\mathbf{785.5}$ & $\mathbf{20.54}$ & $17$ & $724.9$ & $22.34$ & $18$ & $1861.3$ & $24.39$\\
tic-tac-toe & $4$ & $236.7$ & $22.56$ & $\mathbf{19}$ & $\mathbf{124.4}$ & $\mathbf{17.80}$ & $4$ & $181.6$ & $19.40$\\
car\_evaluation & $6$ & $563.7$ & $14.99$ & $\mathbf{19}$ & $\mathbf{348.6}$ & $\mathbf{6.71}$ & $18$ & $733.7$ & $6.19$\\
kr-vs-kp & $7$ & $1475.1$ & $5.79$ & $19$ & $843.8$ & $3.04$ & $\mathbf{20}$ & $\mathbf{913.1}$ & $\mathbf{3.25}$\\
\hline
iris & $55$ & $243.6$ & $0.73$ & $\mathbf{60}$ & $\mathbf{3.7}$ & \textbf{-} & $60$ & $6.8$ & -\\
wine & $48$ & $464.6$ & $1.17$ & $\mathbf{52}$ & $\mathbf{237.0}$ & $\mathbf{0.42}$ & $49$ & $465.2$ & $1.01$\\
plrx & $11$ & $435.5$ & $22.56$ & $\mathbf{36}$ & $\mathbf{102.2}$ & $\mathbf{20.41}$ & $12$ & $384.1$ & $17.74$\\
wpbc & $9$ & $124.9$ & $16.86$ & $\mathbf{34}$ & $\mathbf{90.6}$ & $\mathbf{14.00}$ & $11$ & $296.0$ & $14.40$\\
parkinsons & $24$ & $351.0$ & $5.74$ & $\mathbf{36}$ & $\mathbf{53.3}$ & $\mathbf{4.34}$ & $26$ & $364.4$ & $3.84$\\
sonar & $3$ & $3055.4$ & $15.64$ & $\mathbf{23}$ & $\mathbf{743.6}$ & $\mathbf{12.49}$ & $17$ & $1178.8$ & $12.09$\\
wdbc & $13$ & $298.4$ & $3.95$ & $\mathbf{36}$ & $\mathbf{299.6}$ & $\mathbf{3.91}$ & $18$ & $679.5$ & $3.25$\\
transfusion & $15$ & $281.3$ & $14.15$ & $\mathbf{60}$ & $\mathbf{70.5}$ & \textbf{-} & $60$ & $111.9$ & -\\
banknote & $35$ & $695.9$ & $2.95$ & $\mathbf{60}$ & $\mathbf{30.0}$ & \textbf{-} & $60$ & $54.8$ & -\\
ozone-one & $19$ & $295.9$ & $2.37$ & $\mathbf{39}$ & $\mathbf{787.2}$ & $\mathbf{2.43}$ & $21$ & $236.0$ & $2.42$\\
segmentation & $0$ & - & $50.47$ & $\mathbf{24}$ & $\mathbf{1318.7}$ & $\mathbf{40.33}$ & $0$ & - & $45.95$\\
spambase & $0$ & - & $15.60$ & $\mathbf{20}$ & $\mathbf{1047.5}$ & $\mathbf{12.20}$ & $3$ & $1301.8$ & $13.36$\\
\hline
hepatitis & $33$ & $1145.8$ & $0.78$ & $\mathbf{60}$ & $\mathbf{240.6}$ & \textbf{-} & $43$ & $647.0$ & $0.45$\\
fertility & $31$ & $619.6$ & $2.88$ & $\mathbf{60}$ & $\mathbf{65.6}$ & \textbf{-} & $46$ & $429.5$ & $1.74$\\
ionosphere & $10$ & $584.6$ & $8.36$ & $\mathbf{36}$ & $\mathbf{191.1}$ & $\mathbf{6.78}$ & $11$ & $436.8$ & $7.12$\\
thoracic & $14$ & $156.3$ & $11.69$ & $\mathbf{43}$ & $\mathbf{318.1}$ & $\mathbf{7.04}$ & $21$ & $304.2$ & $5.10$\\
ILPD & $10$ & $419.8$ & $29.31$ & $\mathbf{36}$ & $\mathbf{110.9}$ & $\mathbf{26.83}$ & $12$ & $450.2$ & $23.70$\\
credit & $12$ & $90.0$ & $12.09$ & $\mathbf{36}$ & $\mathbf{186.8}$ & $\mathbf{11.14}$ & $14$ & $495.3$ & $11.13$\\
biodeg & $0$ & - & $24.54$ & $\mathbf{29}$ & $\mathbf{952.0}$ & $\mathbf{19.06}$ & $2$ & $1763.1$ & $21.10$\\
seismic-bumps & $16$ & $9.1$ & $5.55$ & $\mathbf{46}$ & $\mathbf{109.2}$ & $\mathbf{3.01}$ & $31$ & $129.3$ & $2.29$\\
ann-thyroid & $18$ & $357.4$ & $5.66$ & $\mathbf{44}$ & $\mathbf{812.5}$ & $\mathbf{2.34}$ & $29$ & $324.8$ & $2.11$\\
\hline
\end{tabular}
    \caption{Per dataset results. Compares BendOCT with no acceleration techniques, all acceleration techniques, and all acceleration techniques minus path bound cutting planes. Results are aggregated over all available instances, 30 for categorical datasets and 60 for numerical datasets. The solved columns list the number of instances solved to optimality. Time lists the averages solve time in seconds of solved instances, '-' indicates that no instances were solved. Gap lists the average percentage optimality gap in unsolved instances, '-' indicates that all instances were solved. Bold font indicates the best performing variant prioritising the number of instances solve to optimality, the average solve time, and then the average optimality gap.}
    \label{tab:Main Results Table}
\end{table}

\section{Conclusion}\label{sect:Conclusion}
We have proposed acceleration techniques applied for the BendOCT MIP formulation for learning optimal classification trees. This includes a novel derivation of BendOCT by logic-based Benders decomposition which informs a strengthening of the Benders cuts, a solution polishing primal heuristics, a class of valid inequalities exploiting similar points in the training data, and a class of cutting planes which leverage a subroutine for finding small optimal subtrees. We showed experimentally that these extensions provide a 300x speedup, allowing for twice as many instances to be solved to optimality within a one hour time limit. These results demonstrate that ideas from DP can be implemented in a MIP setting using advanced integer programming techniques to significantly accelerate convergence.

There are several directions which can be considered in future work. While the logic underlying the acceleration techniques is applicable to other MIP formulations, whether they can be practically adapted is left for future work. There remain ideas from the DP literature which could be effective in the MIP setting, similar support bounds in particular have found great success and may be critical for scaling MIP formulations to deeper trees. Extending the proposed techniques to other objectives and constraints is important for unlocking the potential of MIP formulations. While the acceleration techniques are powerful, they come at the potential cost of requiring tailoring to suit extensions. The logic of cuts which relate the tree structure to the classification scores would generally need modification to suit alternative objectives. In the case of BendOCT, Benders decomposition provides a flexible framework, with feasibility cuts being a powerful tool for incorporating constraints which cannot be conveniently modelled in the MP. Thus, further work is required to extend the acceleration techniques to more exotic settings.

\section{Acknowledgements}
This research was supported by an Australian Government Research Training Program (RTP) Scholarship which is gratefully acknowledged.

\bibliography{references}

\newpage
\appendix

\section{Appendix}

\subsection{Dataset Statistics}

\begin{table}[h]
    \centering
    \renewcommand*{\arraystretch}{0.85}
    \begin{tabular}{cccccccccc}
\hline
Dataset&$|\mathcal{I}|$&$|\mathcal{F}|$&$|\mathcal{F}^C|$&$|\mathcal{F}^N|$&$|\mathcal{F}_C^B|$&$|\mathcal{F}_{QB-5}^B|$&$|\mathcal{F}_{QT-5}^B|$&$|\mathcal{K}|$&Class Distribution\\
\hline 
soybean-small&47&21&21&0&45&N/A&N/A&4&(21, 36)\\
monk3&122&6&6&0&15&N/A&N/A&2&(49, 51)\\
monk1&124&6&6&0&15&N/A&N/A&2&(50, 50)\\
hayes-roth&132&4&4&0&15&N/A&N/A&3&(23, 39)\\
monk2&169&6&6&0&15&N/A&N/A&2&(38, 62)\\
house-votes-84&232&16&16&0&16&N/A&N/A&2&(47, 53)\\
spect&267&22&22&0&22&N/A&N/A&2&(21, 79)\\
breast-cancer&277&9&9&0&38&N/A&N/A&2&(29, 71)\\
balance-scale&625&4&4&0&20&N/A&N/A&3&(8, 46)\\
tic-tac-toe&958&9&9&0&27&N/A&N/A&2&(35, 65)\\
car-evaluation&1728&6&6&0&19&N/A&N/A&4&(4, 70)\\
kr-vs-kp&3196&36&36&0&38&N/A&N/A&2&(48, 52)\\
\hline 
iris&150&4&0&4&N/A&20&16&3&(33, 33)\\
wine&178&13&0&13&N/A&65&52&3&(27, 40)\\
plrx&182&12&0&12&N/A&60&48&2&(29, 71)\\
wpbc&194&32&0&32&N/A&158&127&2&(24, 76)\\
parkinsons&195&22&0&22&N/A&110&88&2&(25, 75)\\
sonar&208&60&0&60&N/A&300&240&2&(47, 53)\\
wdbc&569&30&0&30&N/A&150&120&2&(37, 63)\\
transfusion&748&4&0&4&N/A&18&16&2&(24, 76)\\
banknote&1372&4&0&4&N/A&20&16&2&(44, 56)\\
ozone-one&1848&72&0&72&N/A&357&286&2&(3, 97)\\
segmentation&2310&18&0&18&N/A&82&66&7&(14, 14)\\
spambase&4601&57&0&57&N/A&92&89&2&(39, 61)\\
\hline 
hepatitis&80&19&13&6&N/A&43&37&2&(16, 84)\\
fertility&100&9&6&3&N/A&27&25&2&(12, 88)\\
ionosphere&351&33&1&32&N/A&157&129&2&(36, 64)\\
thoracic&470&16&13&3&N/A&39&36&2&(15, 85)\\
ILPD&579&10&1&9&N/A&45&36&2&(28, 72)\\
credit&653&15&10&5&N/A&82&78&2&(45, 55)\\
biodeg&1055&41&9&32&N/A&143&127&2&(34, 66)\\
seismic-bumps&2584&15&5&10&N/A&43&39&2&(7, 93)\\
ann-thyroid&7200&21&15&6&N/A&45&39&3&(2, 93)\\
\hline
\end{tabular}

    \caption{Table of datasets used in computational experiments with basic statistics. Datasets are separated into those with purely categorical, purely continuous, and mixed categorical and continuous features. For each dataset the number of categorical and continuous features are listed along with the number of binary features under each binary encoding (one-hot, quantile bucketing, or quantile thresholding). The class distribution column lists the proportions of the minority and majority classes.}
    \label{tab: Dataset Statistics}
\end{table}

\newpage

\section{Online Appendix}

\subsection{Optimal Depth 2 Tree Subroutine}\label{Appendix:Optimal Depth 2 Tree Subroutine}
MuRTree is a dynamic programming formulation for learning optimal classification trees \citep{demirovicMurTreeOptimalDecision2022}. A key element is a subroutine for efficiently finding optimal depth 2 trees. It is presented for the case of depth two binary classification trees, but can be trivially extended for multi-class classification. In this Appendix we detail an extension to find an optimal sparse subtree given a per leaf penalty on classification score $\bar{\lambda} = \lambda |\mathcal{I}|$. The original subroutine implicitly calculates all the information required to find the optimal sparse subtree, we only modify it to explicitly find this sparse subtree. The key ideas are as follows:
\begin{itemize}
    \item It is possible to efficiently calculate frequency counters for how often samples of each class fall into each node in the subtree for each combination of parent and left/right child features.
    \item Given these frequency counters, it is possible to quickly calculate the classification score in each subtree for each possible combination of branch decisions for each valid tree topology.
    \item The optimal subtree for each number of leaf nodes can be efficiently calculated, after which the optimal sparse subtree is whichever one maximises the classification score subject to leaf penalty $\bar{\lambda}$.
\end{itemize}

\begin{figure}[htb!]
    \centering
    \includegraphics[width=0.5\linewidth]{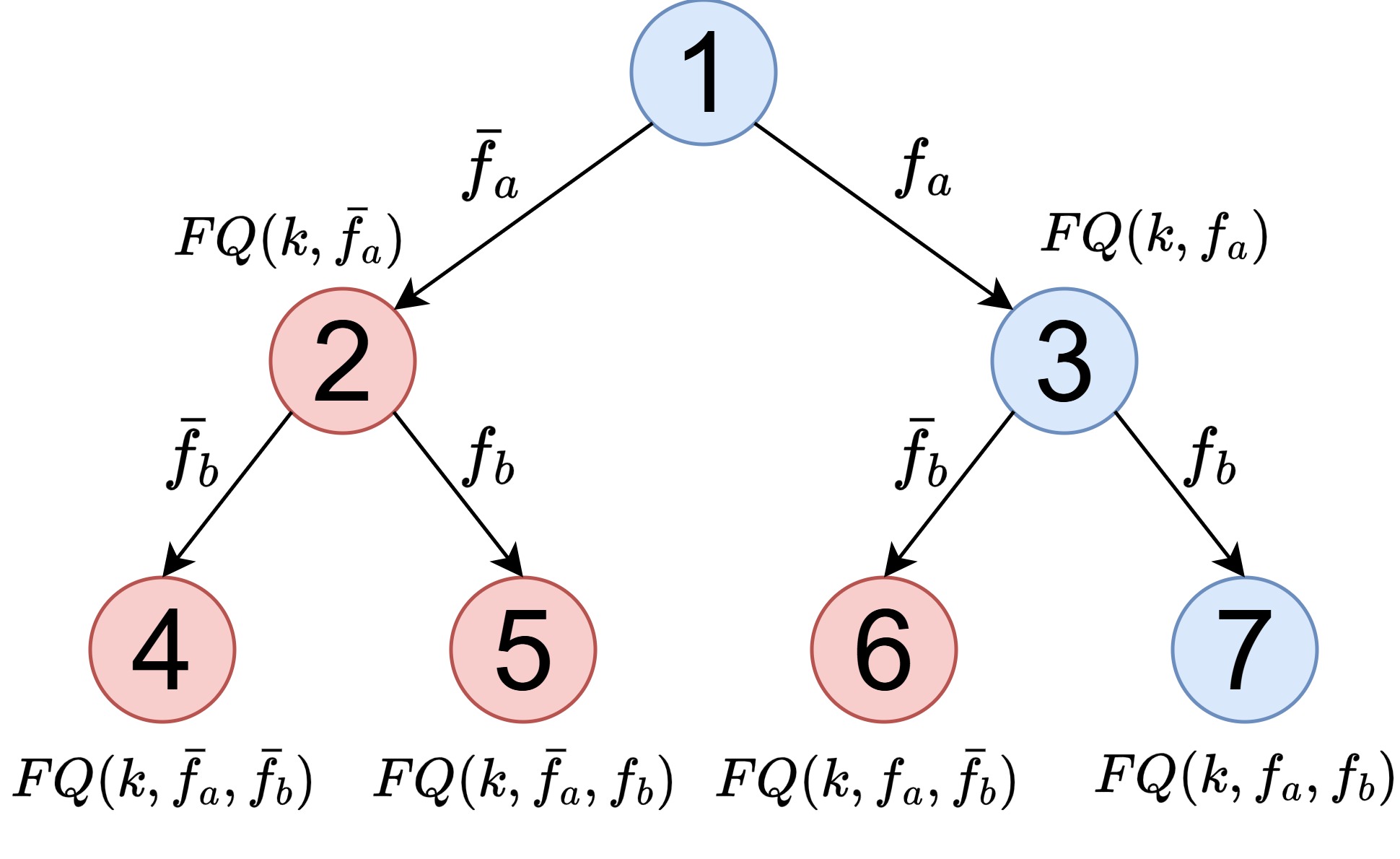}
    \caption{Subroutine Frequency Counter Construction. Blue nodes indicate counters computed by brute force, red nodes indicate derived counters.}
    \label{fig:subroutine diagram}
\end{figure}

To calculate the frequency counters, denote the set of samples in class $k$ as $\mathcal{D}^k$ such that all data in the root $\mathcal{D} = \cup_{k \in K} \mathcal{D}^k$. Denote the number of samples with class $k$ that flow right from the parent node branching on feature $f_a \in F$ as $FQ(k,f_a) = |\{i \in D^k :  x_{f_a}^i = 1\}|$. Similarly denote $FQ(k,f_a, f_b) = |\{i \in D^k :  x_{f_a}^i = 1 \; and \; x_{f_b}^i = 1\}|$ as the number of samples branching right twice, first on $f_a$ then on $f_b$. We use the notation $\bar{f}$ to denote that feature $f$ is branched left on, E.g. $FQ(k,\bar{f}_a, f_b) = |\{i \in D^k :  x_{f_a}^i = 0 \; and \; x_{f_b}^i = 1\}|$ is the number of samples which branch left on feature $f_a$ then right on feature $f_b$. Given $FQ(k,f_a)$ and $FQ(k,f_a, f_b)$, it is possible to quickly compute $FQ(k,\bar{f}_a)$, $FQ(k,\bar{f}_a, f_b)$, $FQ(k,\bar{f}_a, \bar{f}_b)$, and $FQ(k,f_a, \bar{f}_b)$. Each of these frequency counters corresponds to either one of the children of the root node, or one of the terminal nodes of the subtree as in Figure \ref{fig:subroutine diagram}. Given $FQ(k,f_a)$ and $FQ(k,f_a,f_b)$ the other frequency counters are computed as follows:

\begin{align}
    FQ(k, \bar{f}_a) &= |\mathcal{D}^k| - FQ(k,f_a) \\ \nonumber \\
    FQ(k, f_a, \bar{f}_b) &= FQ(k, f_a) - FQ(k,f_a, f_b) \\ \nonumber  \\
    FQ(k, \bar{f}_a, f_b) &= FQ(k, f_b) - FQ(k,f_a, f_b) \\ \nonumber  \\
    FQ(k, \bar{f}_a, \bar{f}_b) &= FQ(k, \bar{f}_a) - FQ(k, \bar{f}_a, f_b)
\end{align}

We use these frequency counters to calculate the number of correctly classified points in each potential leaf node given any possible combination of branch decisions in ancestor nodes. We denote the classification score in the root node by $CS$, in the children of the root node given branch decision $f_a$ by $CS(f_a)$, and in the terminal nodes given branch decision $f_a$ and $f_b$ by $CS(f_a,f_b)$. For each combination of $(f_a, \bar{f}_a, f_b, \bar{f}_b)$, the classification scores are:


\begin{align}
    CS &= \max_{k \in K} |D^k| \\
    CS(f_a) &= \max_{k \in K} FQ(k,f_a) \\
    CS(f_a,f_b) &= \max_{k \in K} FQ(k,f_a,f_b) 
\end{align}

Given these classification scores we find the optimal subtree for each number of leaf nodes. The optimal single node subtree obviously has a classification score of $CS_1= CS$. The optimal two node subtree has a classification score of:

\begin{equation}
    CS_2 = \max_{f_a \in F}CS(f_a) + CS(\bar{f}_a)
\end{equation}

There are two possible three node subtrees differing in which child node is branched on. We first calculate the classification scores in the fully grown left and right subtrees over all possible combinations of left and right child features:

\begin{align}
    CS^L(f_a) &= \max_{f_b \in F} CS(\bar{f}_a,f_b) + CS(\bar{f}_a,\bar{f_b}) \\
    CS^R(f_a) &= \max_{f_b \in F} CS(f_a,f_b) + CS(f_a,\bar{f_b})
\end{align}

Given this we can find the optimal three leaf subtrees given that the left or right child nodes are fully expanded. The optimal three leaf subtree is then whichever of the two has a higher classification score.

\begin{align}
    CS_3^L &= \max_{f_a \in F} CS(f_a) + CS^L(f_a) \\
    CS_3^R &= \max_{f_a \in F} CS(\bar{f}_a) + CS^R(f_a) \\
    CS_3 &= \max \{ CS_3^L , CS_3^R \}
\end{align}

The classification score in the optimal four leaf subtree is calculated in a similar fashion:

\begin{equation*}
    CS_4 = \max_{f_a \in F} CS^L(f_a) + CS^R(f_a)
\end{equation*}

Finally, the number of correctly classified points in the optimal subtree is:

\begin{equation}
    \max_{H \in  \{1,2,3,4\}}
    CS_H - \bar{\lambda} H
\end{equation}

Critically the subroutine never optimises over combinations branch decisions at the three branch nodes. Smaller optimisations are done separately and combined to construct optimal subtrees for each number of leaf nodes. The time complexity is $\mathcal{O} (|\mathcal{I}|\cdot |F|^2)$ with $\mathcal{O} (|F|^2)$ space used. We note that there may be multiple subtrees with the same optimal objective, which theoretically could differ in the number of leaves. In these situations our implementation selects the lexicographically smallest subtree in the branch features, but it is trivial to modify the subroutine to select for specific qualities within the set of symmetric solutions.

We implement the computationally expensive operations of calculating the frequency counters and maximising over frequency counters and classification scores using Numpy arrays. We found that in most cases the time spent in the subroutine was negligible compared to the overall optimisation time, and as such did not pursue further optimisations. \citet{demirovicMurTreeOptimalDecision2022} describe two further optimisations. The first is to note that binary feature vectors tend to be sparse, which can be exploited by representing the feature vectors as sparse vectors and iterating over the sparse entries to update the frequency counters. The second is to store old frequency counters which can be incrementally updated to calculate optimal subtrees for datasets which differ only in a small number of samples. This is incredibly effective in their dynamic programming setting since the subroutine is often called to optimise subtrees which differ only in the branch feature in the parent node of the subtree, which will often result in sequentially optimising over similar subsets of the samples. This is less likely to occur in the MIP setting, as such benefits are unlikely to justify the computation and memory overhead for storing and incrementally updating the frequency counters.

\subsection{Extended Results}\label{Appendix: Extended Results}
For some of the acceleration techniques proposed in this work we have presented multiple variants. Section \ref{sect: Main Results} presents results for the best performing variants. In this Appendix we justify the use of specific cut configurations by computational experiments, and present per dataset results for each acceleration technique. All experiments follow the experimental setup described in Section \ref{sect: Computational Experiments}. All figures should be interpreted as follows: the left panel indicates how many instances each model solved to optimality within the given time, the right panel indicates for how many instances each model found a solution within the given optimality gap. All tables report per dataset results aggregated over all available instances, meaning all possible combinations of maximum tree depth, $\lambda$, and binary encoding for each dataset.

\subsection{Strengthened Benders Cuts}
We report per dataset results for the strengthened Benders cuts in Table \ref{tab:Aggregated Strengthened Cuts Results}. The strengthened cuts provide a consistent marginal benefit for most datasets, except for the sonar, banknote, and hepatitis datasets for which the strengthened cuts allow for an extra 14, 14, and 9 instances to be solved respectively.

\begin{table}[htb!]
    \centering
    \renewcommand*{\arraystretch}{0.85}
    \begin{tabular}{ccccccc}
\hline
\multirow{2}{*}{Dataset} & \multicolumn{3}{c}{BendOCT} & \multicolumn{3}{c}{Strengthened Benders Cuts}\\
  & Solved & Time & Gap & Solved & Time & Gap\\
\hline
soybean-small & $30$ & $0.2$ & - & $\mathbf{30}$ & $\mathbf{0.2}$ & \textbf{-}\\
monk3 & $25$ & $258.5$ & $0.30$ & $\mathbf{27}$ & $\mathbf{517.8}$ & $\mathbf{0.78}$\\
monk1 & $\mathbf{30}$ & $\mathbf{6.4}$ & \textbf{-} & $30$ & $6.5$ & -\\
hayes-roth & $\mathbf{20}$ & $\mathbf{353.6}$ & $\mathbf{2.72}$ & $18$ & $132.1$ & $4.72$\\
monk2 & $15$ & $1642.6$ & $17.89$ & $\mathbf{18}$ & $\mathbf{1003.7}$ & $\mathbf{20.04}$\\
house-votes-84 & $26$ & $257.3$ & $0.27$ & $\mathbf{27}$ & $\mathbf{348.0}$ & $\mathbf{0.23}$\\
spect & $18$ & $956.7$ & $3.63$ & $\mathbf{20}$ & $\mathbf{568.3}$ & $\mathbf{4.23}$\\
breast-cancer & $7$ & $468.7$ & $16.26$ & $\mathbf{7}$ & $\mathbf{337.8}$ & $\mathbf{13.43}$\\
balance-scale & $18$ & $785.5$ & $20.54$ & $\mathbf{18}$ & $\mathbf{696.6}$ & $\mathbf{14.69}$\\
tic-tac-toe & $4$ & $236.7$ & $22.56$ & $\mathbf{4}$ & $\mathbf{165.4}$ & $\mathbf{19.99}$\\
car\_evaluation & $\mathbf{6}$ & $\mathbf{563.7}$ & $\mathbf{14.99}$ & $6$ & $581.0$ & $14.34$\\
kr-vs-kp & $\mathbf{7}$ & $\mathbf{1475.1}$ & $\mathbf{5.79}$ & $4$ & $1483.1$ & $6.74$\\
\hline
iris & $55$ & $243.6$ & $0.73$ & $\mathbf{58}$ & $\mathbf{301.4}$ & $\mathbf{0.58}$\\
wine & $\mathbf{48}$ & $\mathbf{464.6}$ & $\mathbf{1.17}$ & $48$ & $470.0$ & $1.03$\\
plrx & $11$ & $435.5$ & $22.56$ & $\mathbf{12}$ & $\mathbf{368.9}$ & $\mathbf{18.85}$\\
wpbc & $9$ & $124.9$ & $16.86$ & $\mathbf{11}$ & $\mathbf{304.6}$ & $\mathbf{15.27}$\\
parkinsons & $24$ & $351.0$ & $5.74$ & $\mathbf{25}$ & $\mathbf{243.3}$ & $\mathbf{4.60}$\\
sonar & $3$ & $3055.4$ & $15.64$ & $\mathbf{17}$ & $\mathbf{1317.9}$ & $\mathbf{12.89}$\\
wdbc & $13$ & $298.4$ & $3.95$ & $\mathbf{17}$ & $\mathbf{556.0}$ & $\mathbf{3.47}$\\
transfusion & $15$ & $281.3$ & $14.15$ & $\mathbf{18}$ & $\mathbf{679.3}$ & $\mathbf{13.67}$\\
banknote & $35$ & $695.9$ & $2.95$ & $\mathbf{49}$ & $\mathbf{797.7}$ & $\mathbf{4.16}$\\
ozone-one & $19$ & $295.9$ & $2.37$ & $\mathbf{21}$ & $\mathbf{216.8}$ & $\mathbf{2.45}$\\
segmentation & $0$ & - & $50.47$ & $\mathbf{0}$ & \textbf{-} & $\mathbf{49.94}$\\
spambase & $0$ & - & $15.60$ & $\mathbf{1}$ & $\mathbf{1096.4}$ & $\mathbf{15.36}$\\
\hline
hepatitis & $33$ & $1145.8$ & $0.78$ & $\mathbf{42}$ & $\mathbf{587.3}$ & $\mathbf{0.51}$\\
fertility & $31$ & $619.6$ & $2.88$ & $\mathbf{38}$ & $\mathbf{243.6}$ & $\mathbf{2.63}$\\
ionosphere & $10$ & $584.6$ & $8.36$ & $\mathbf{11}$ & $\mathbf{441.7}$ & $\mathbf{7.82}$\\
thoracic & $14$ & $156.3$ & $11.69$ & $\mathbf{16}$ & $\mathbf{437.9}$ & $\mathbf{10.29}$\\
ILPD & $10$ & $419.8$ & $29.31$ & $\mathbf{10}$ & $\mathbf{310.4}$ & $\mathbf{27.30}$\\
credit & $12$ & $90.0$ & $12.09$ & $\mathbf{13}$ & $\mathbf{380.6}$ & $\mathbf{11.93}$\\
biodeg & $0$ & - & $24.54$ & $\mathbf{0}$ & \textbf{-} & $\mathbf{23.54}$\\
seismic-bumps & $\mathbf{16}$ & $\mathbf{9.1}$ & $\mathbf{5.55}$ & $16$ & $11.5$ & $5.55$\\
ann-thyroid & $\mathbf{18}$ & $\mathbf{357.4}$ & $\mathbf{5.66}$ & $17$ & $276.9$ & $5.93$\\
\hline
\end{tabular}
    \caption{Strengthened Benders cuts per dataset results. The solved columns list the number of instances solved to optimality. Time lists the averages solve time in seconds of solved instances, ’-’ indicates that no instances were solved. Gap lists the average percentage optimality gap in unsolved instances, ’-’ indicates that all instances were solved. Bold font indicates the best performing variant prioritising the number of instances solve to optimality, the average solve time, followed by the average optimality gap.}
    \label{tab:Aggregated Strengthened Cuts Results}
\end{table}

\subsection{Solution Polishing}
We report per dataset results for the solution polishing primal heuristic in Table \ref{tab:Aggregated Solution Polishing Results}. The primal heuristic rarely results in more instances being solved to optimality, but in most cases it slightly increases the average solve time while marginally improving the optimality gap. This indicates that the primal heuristic does often find slightly improved solutions, but since it does little to improve the upper bound it rarely assists in proving optimality and slightly increases the solve time in instances which would have been solved to optimality regardless.

\begin{table}[htb!]
    \centering
    \renewcommand*{\arraystretch}{0.85}
    \begin{tabular}{ccccccc}
\hline
\multirow{2}{*}{Dataset} & \multicolumn{3}{c}{BendOCT} & \multicolumn{3}{c}{Solution Polshing}\\
  & Solved & Time & Gap & Solved & Time & Gap\\
\hline
soybean-small & $\mathbf{30}$ & $\mathbf{0.2}$ & \textbf{-} & $30$ & $0.5$ & -\\
monk3 & $\mathbf{25}$ & $\mathbf{258.5}$ & $\mathbf{0.30}$ & $25$ & $302.3$ & $0.31$\\
monk1 & $30$ & $6.4$ & - & $\mathbf{30}$ & $\mathbf{5.6}$ & \textbf{-}\\
hayes-roth & $20$ & $353.6$ & $2.72$ & $\mathbf{21}$ & $\mathbf{549.0}$ & $\mathbf{3.11}$\\
monk2 & $15$ & $1642.6$ & $17.89$ & $\mathbf{16}$ & $\mathbf{1784.1}$ & $\mathbf{17.31}$\\
house-votes-84 & $26$ & $257.3$ & $0.27$ & $\mathbf{26}$ & $\mathbf{217.1}$ & $\mathbf{0.27}$\\
spect & $18$ & $956.7$ & $3.63$ & $\mathbf{18}$ & $\mathbf{898.8}$ & $\mathbf{3.56}$\\
breast-cancer & $7$ & $468.7$ & $16.26$ & $\mathbf{7}$ & $\mathbf{465.0}$ & $\mathbf{15.87}$\\
balance-scale & $\mathbf{18}$ & $\mathbf{785.5}$ & $\mathbf{20.54}$ & $18$ & $791.4$ & $19.63$\\
tic-tac-toe & $\mathbf{4}$ & $\mathbf{236.7}$ & $\mathbf{22.56}$ & $4$ & $241.4$ & $21.56$\\
car\_evaluation & $\mathbf{6}$ & $\mathbf{563.7}$ & $\mathbf{14.99}$ & $6$ & $565.2$ & $15.02$\\
kr-vs-kp & $\mathbf{7}$ & $\mathbf{1475.1}$ & $\mathbf{5.79}$ & $7$ & $1485.6$ & $5.13$\\
\hline
iris & $\mathbf{55}$ & $\mathbf{243.6}$ & $\mathbf{0.73}$ & $53$ & $130.2$ & $0.58$\\
wine & $\mathbf{48}$ & $\mathbf{464.6}$ & $\mathbf{1.17}$ & $47$ & $334.5$ & $0.90$\\
plrx & $\mathbf{11}$ & $\mathbf{435.5}$ & $\mathbf{22.56}$ & $11$ & $440.6$ & $20.94$\\
wpbc & $9$ & $124.9$ & $16.86$ & $\mathbf{9}$ & $\mathbf{123.0}$ & $\mathbf{15.77}$\\
parkinsons & $24$ & $351.0$ & $5.74$ & $\mathbf{24}$ & $\mathbf{346.9}$ & $\mathbf{4.62}$\\
sonar & $3$ & $3055.4$ & $15.64$ & $\mathbf{3}$ & $\mathbf{3025.9}$ & $\mathbf{14.52}$\\
wdbc & $\mathbf{13}$ & $\mathbf{298.4}$ & $\mathbf{3.95}$ & $13$ & $339.4$ & $3.61$\\
transfusion & $\mathbf{15}$ & $\mathbf{281.3}$ & $\mathbf{14.15}$ & $15$ & $287.3$ & $14.16$\\
banknote & $35$ & $695.9$ & $2.95$ & $\mathbf{39}$ & $\mathbf{792.9}$ & $\mathbf{3.05}$\\
ozone-one & $\mathbf{19}$ & $\mathbf{295.9}$ & $\mathbf{2.37}$ & $19$ & $308.2$ & $2.33$\\
segmentation & $0$ & - & $50.47$ & $\mathbf{0}$ & \textbf{-} & $\mathbf{49.03}$\\
spambase & $0$ & - & $15.60$ & $\mathbf{0}$ & \textbf{-} & $\mathbf{15.10}$\\
\hline
hepatitis & $33$ & $1145.8$ & $0.78$ & $\mathbf{36}$ & $\mathbf{1248.4}$ & $\mathbf{0.82}$\\
fertility & $31$ & $619.6$ & $2.88$ & $\mathbf{34}$ & $\mathbf{764.1}$ & $\mathbf{2.52}$\\
ionosphere & $\mathbf{10}$ & $\mathbf{584.6}$ & $\mathbf{8.36}$ & $10$ & $595.3$ & $7.55$\\
thoracic & $\mathbf{14}$ & $\mathbf{156.3}$ & $\mathbf{11.69}$ & $14$ & $158.9$ & $11.42$\\
ILPD & $\mathbf{10}$ & $\mathbf{419.8}$ & $\mathbf{29.31}$ & $10$ & $455.7$ & $27.64$\\
credit & $\mathbf{12}$ & $\mathbf{90.0}$ & $\mathbf{12.09}$ & $12$ & $90.4$ & $11.89$\\
biodeg & $0$ & - & $24.54$ & $\mathbf{0}$ & \textbf{-} & $\mathbf{23.05}$\\
seismic-bumps & $\mathbf{16}$ & $\mathbf{9.1}$ & $\mathbf{5.55}$ & $16$ & $9.3$ & $5.56$\\
ann-thyroid & $\mathbf{18}$ & $\mathbf{357.4}$ & $\mathbf{5.66}$ & $18$ & $358.5$ & $5.60$\\
\hline
\end{tabular}
    \caption{Solution polishing primal heuristic per dataset results. The solved columns list the number of instances solved to optimality. Time lists the averages solve time in seconds of solved instances, ’-’ indicates that no instances were solved. Gap lists the average percentage optimality gap in unsolved instances, ’-’ indicates that all instances were solved. Bold font indicates the best performing variant prioritising the number of instances solve to optimality, the average solve time, followed by the average optimality gap.}
    \label{tab:Aggregated Solution Polishing Results}
\end{table}

\subsubsection{Equivalent Point Cuts}
The equivalent point inequalities introduced in Section \ref{sect: Equivalent Point Valid Inequalities} permit a number of variants. We present results over combinations of the following:

\begin{itemize}
    \item \textbf{Variant of }$\mathcal{H}$ - We test the basic (Ba), chain (Ch), and recursive (Re) variants.
    \item \textbf{Variant of }$\mathcal{G}$ - By default we use the basic version which naively bounds classification scores. We also test the group selection (GS) variant.
    \item \textbf{Constraint Disaggregation} - $\mathcal{H}^{Chain}$ and $\mathcal{H}^{Recursive}$ both have constraints which can be disaggregated allowing for the integrality of associated variables to relaxed. We assess the impact of this disaggregation (DA).
    \item \textbf{Maximum cardinality of }$F^*$ - Adding all valid EQP cuts for a given dataset overburdens the model and is counter-productive. We limit the number of EQP cuts by only included EQP sets such that $|F^*| \leq FR$ for some upper bound $FR$. We test upper bounds of $0$, $1$, and $2$.
\end{itemize}

\begin{figure}[htb!]
    \centering
    \begin{subfigure}{0.4\textwidth}
        \centering
         \includegraphics[width=\textwidth]{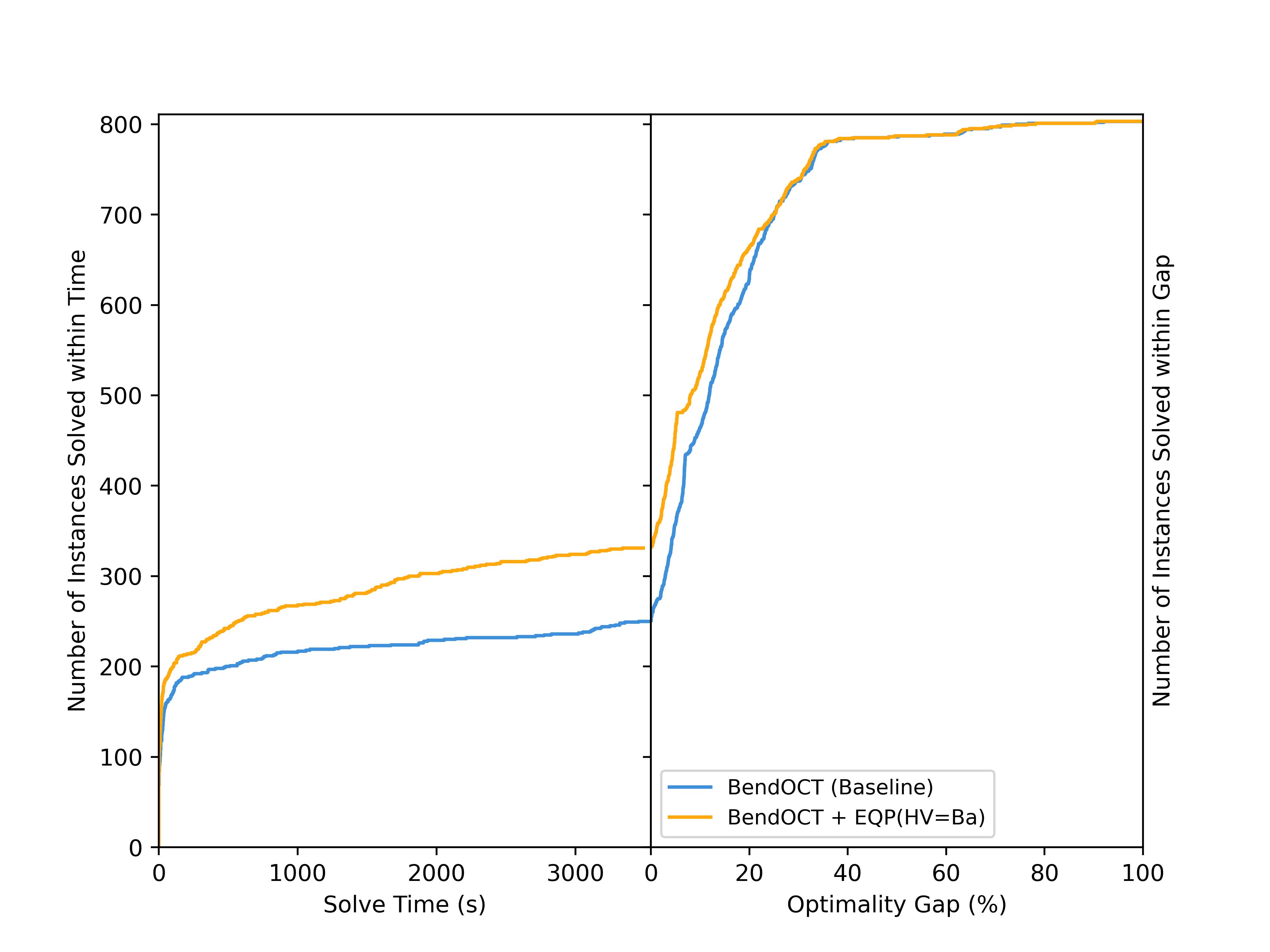}
         \caption{$|F^*| = 0$}
    \end{subfigure}
    \begin{subfigure}{0.4\textwidth}
        \centering
         \includegraphics[width=\textwidth]{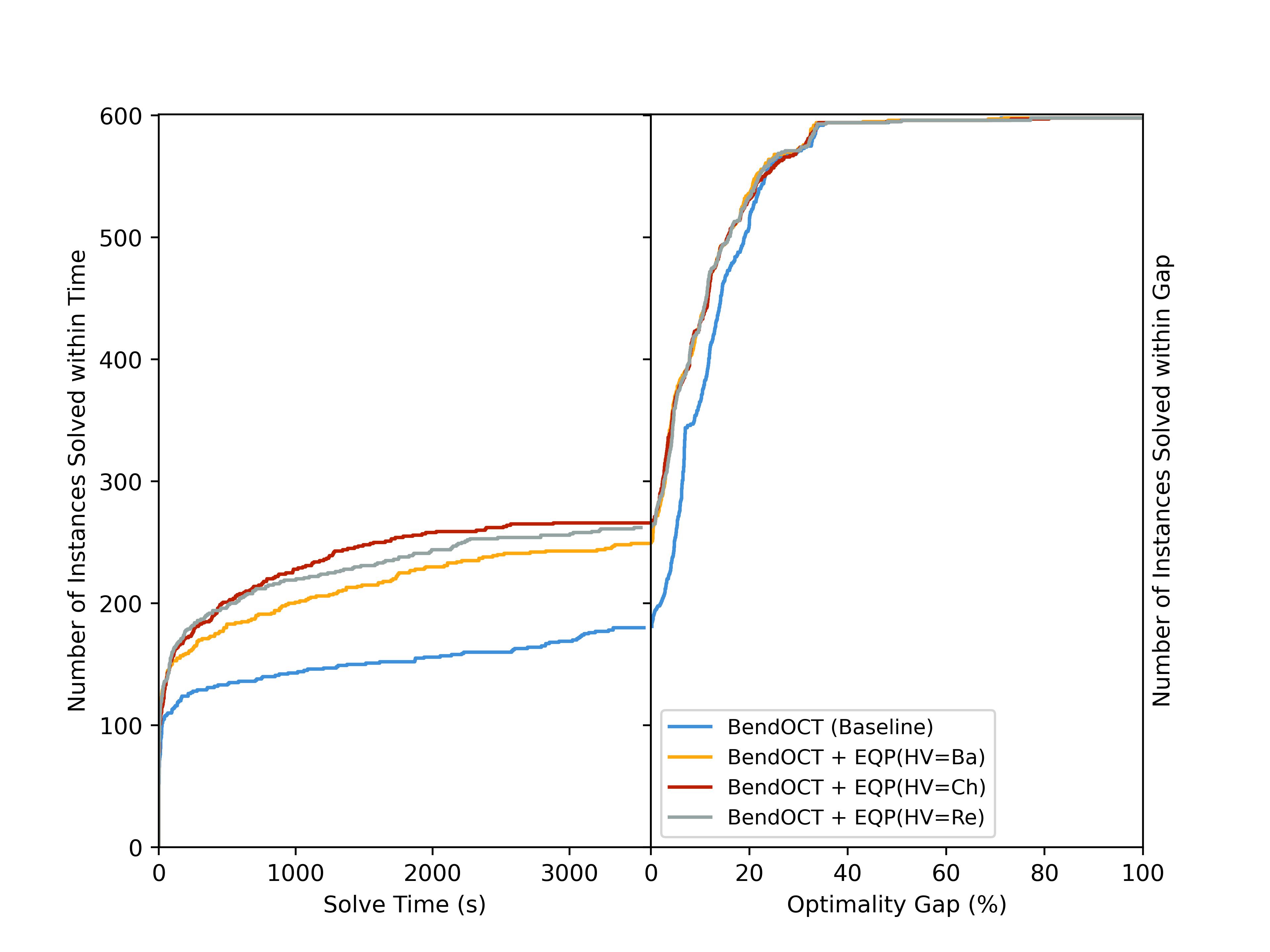}
         \caption{$|F^*| \leq 1$}
    \end{subfigure}
    \hfill
    \begin{subfigure}{0.4\textwidth}
         \centering
         \includegraphics[width=\textwidth]{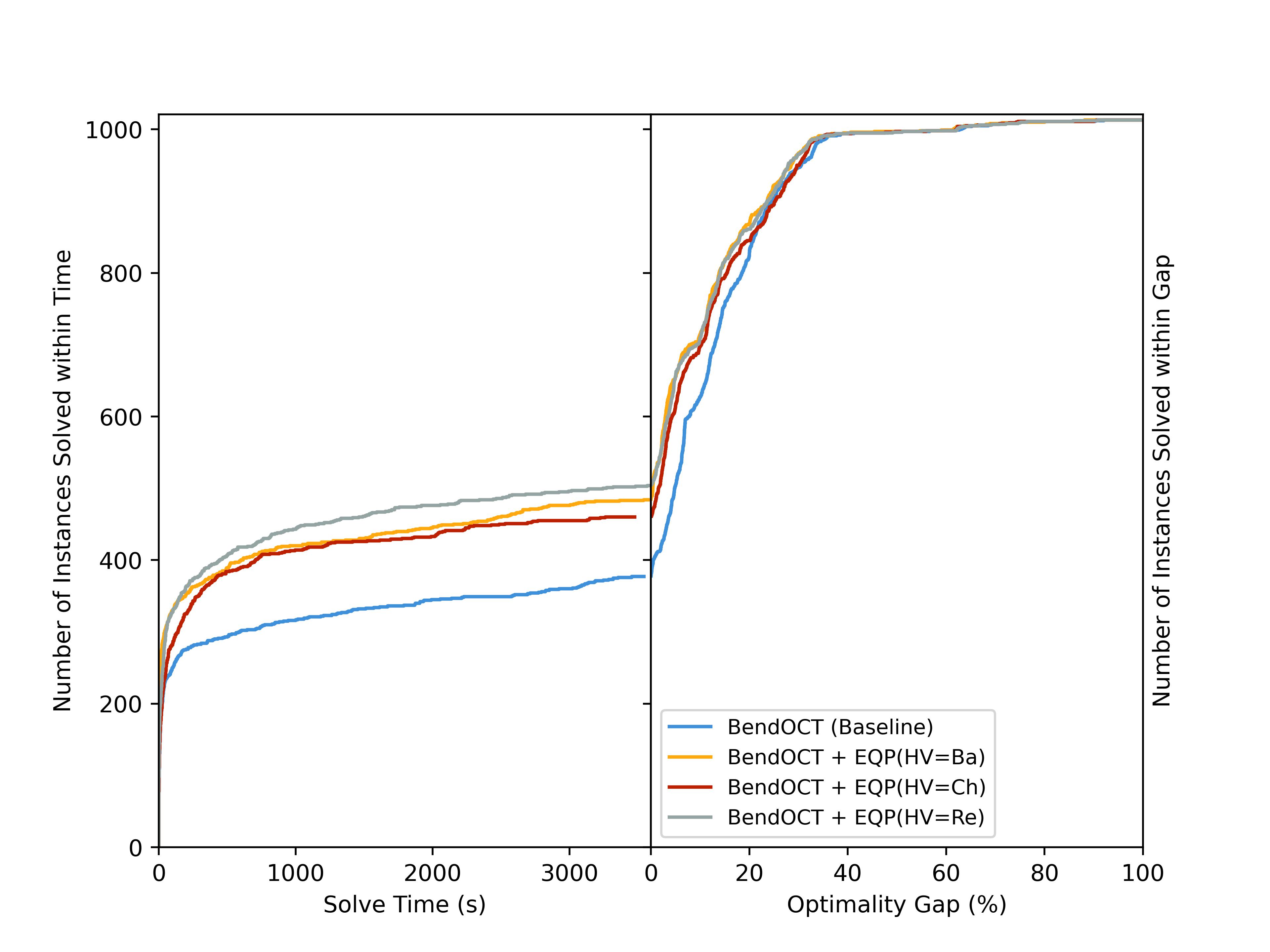}
         \caption{$|F^*| \leq 2$}
    \end{subfigure}
    \caption{Equivalent point inequalities computational results for variants of $\mathcal{H}$ (basic, chain, or recursive) without constraint disaggregation or group selection constraints.}
    \label{fig:EQP H Variant Baseline}
\end{figure}

Figure \ref{fig:EQP H Variant Baseline} presents a basic baseline comparing the different variants of $\mathcal{H}$ without group selection constraints or constraint disaggregation. Results are split into three subplots according to the maximum cardinality of $|F^*|$ allowed when generating EQP sets. We exclude datasets entirely which do not have EQP sets when $|F^*| \leq 2$, for each level of $FR \in \{0,1,2\}$ we only run experiments on instances for which $FR$ contains more EQP sets than for $FR-1$. If an instance is excluded it means that it would yield the same results as using a smaller $FR$ or the baseline BendOCT model. We find $\mathcal{H}^{Basic}$ to be more effective for easier instances which can be solved to optimality within 100 seconds, for harder instances $\mathcal{H}^{Recursive}$ is most effective. There is less variation between the three for instances which cannot be solved to optimality. All three variants dominate vanilla BendOCT.

\begin{figure}[htb!]
    \centering
    \begin{subfigure}[t]{0.4\textwidth}
        \centering
        \includegraphics[width=\textwidth]{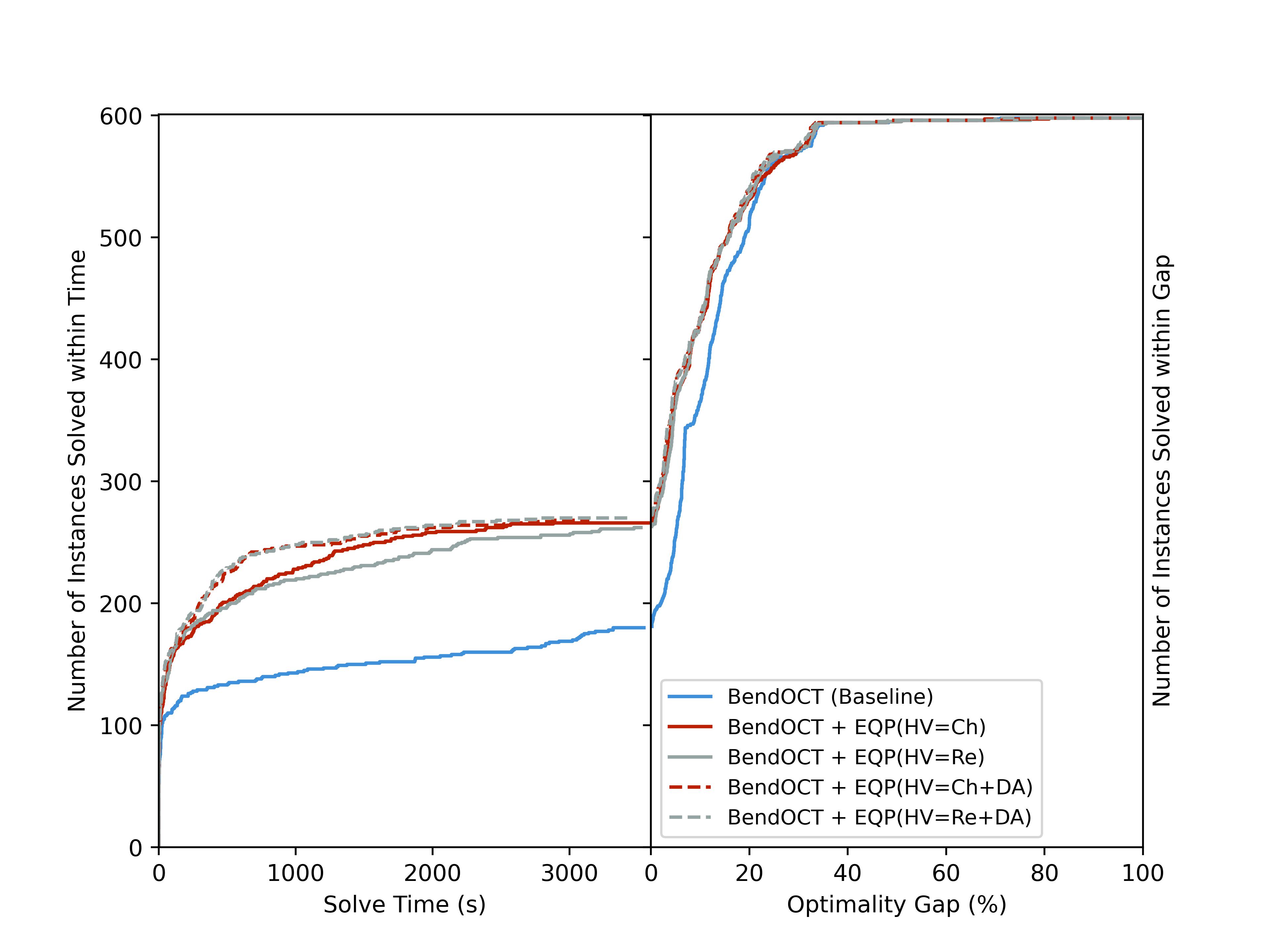}
        \caption{$|F^*| \leq 1$}
    \end{subfigure}
    \begin{subfigure}[t]{0.4\textwidth}
        \centering
        \includegraphics[width=\textwidth]{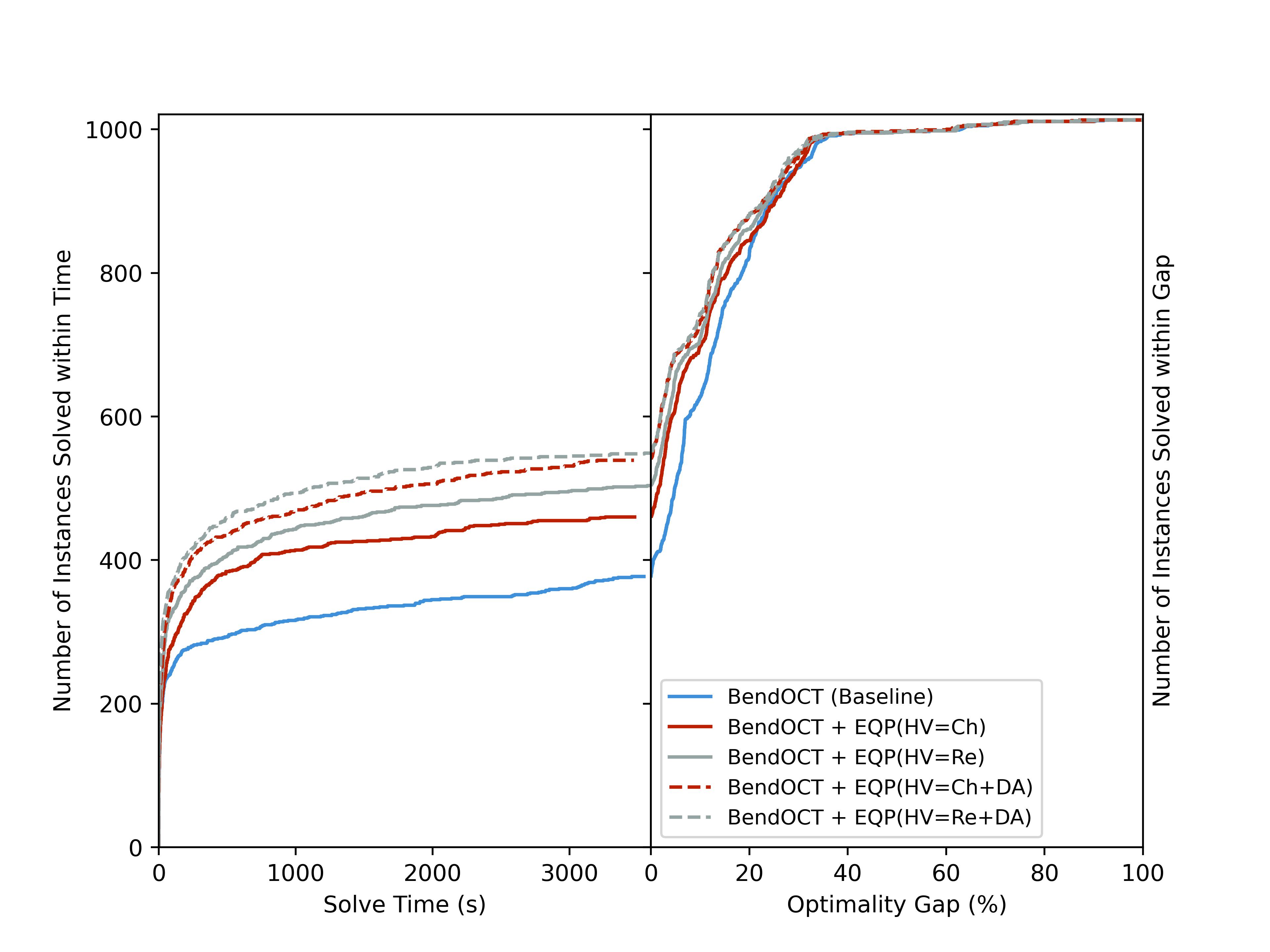}
        \caption{$|F^*| \leq 2$}
    \end{subfigure}
    \caption{Equivalent point inequalities computational results assessing the effect of constraint disaggregation (DA) for chain and recursive variants of $\mathcal{H}$. Disaggregated variants are dashed.}
    \label{fig:EQP Disaggregation Baseline}
\end{figure}

\begin{figure}[htb!]
    \centering
    \begin{subfigure}{0.4\textwidth}
        \centering
         \includegraphics[width=\textwidth]{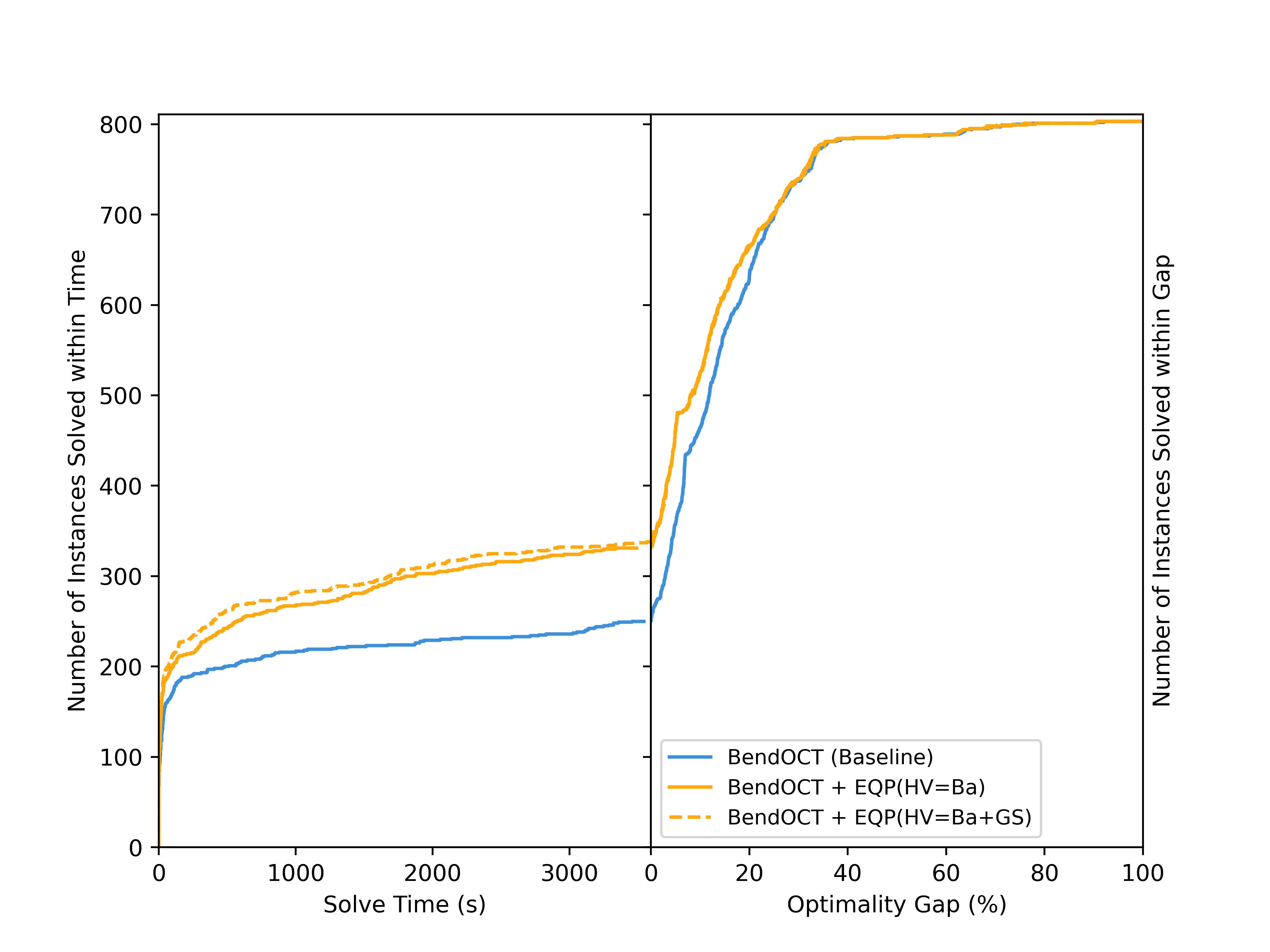}
         \caption{$|F^*| = 0$}
    \end{subfigure}
    \begin{subfigure}{0.4\textwidth}
        \centering
         \includegraphics[width=\textwidth]{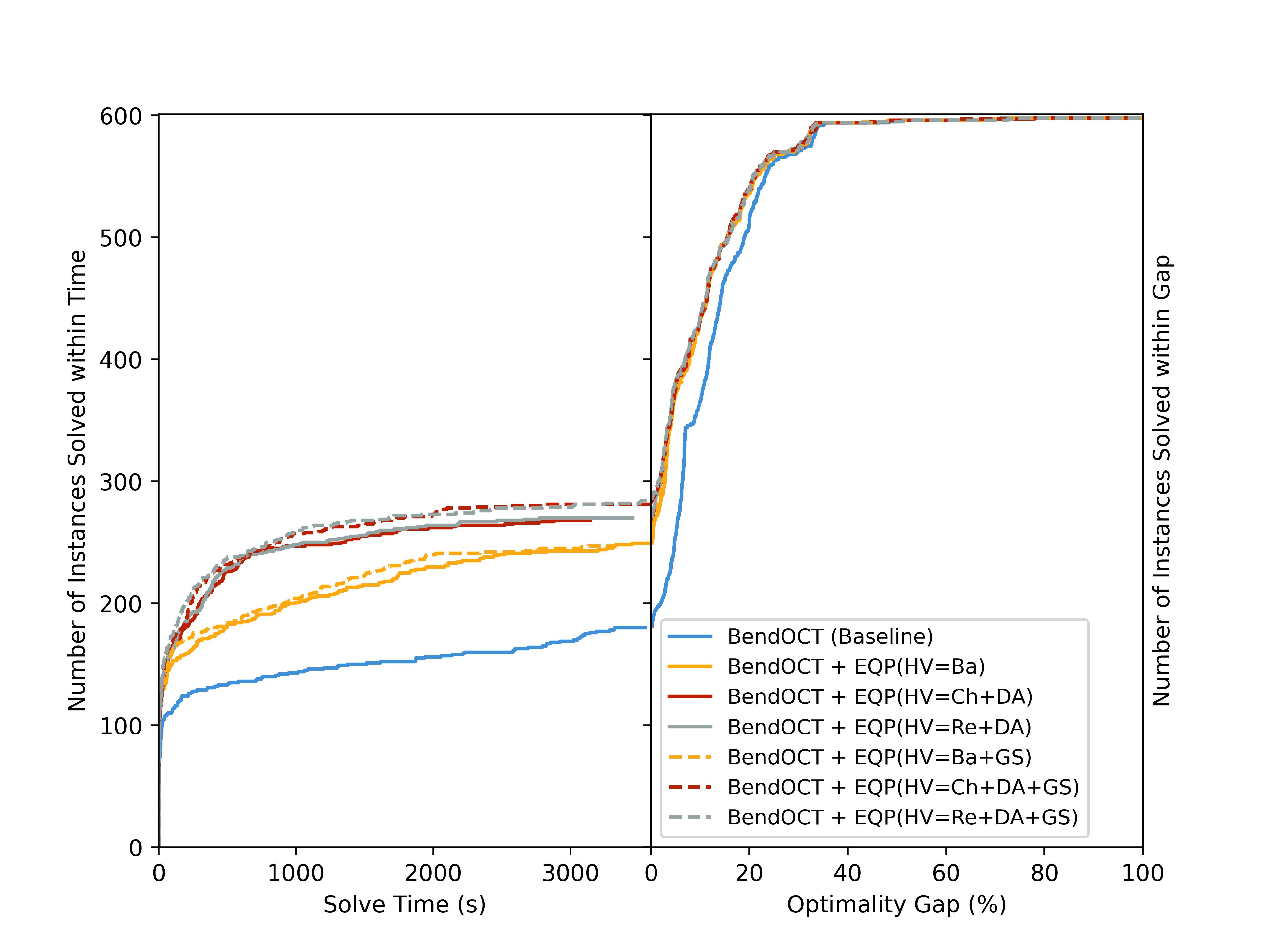}
         \caption{$|F^*| \leq 1$}
    \end{subfigure}
    \hfill
    \begin{subfigure}{0.4\textwidth}
         \centering
         \includegraphics[width=\textwidth]{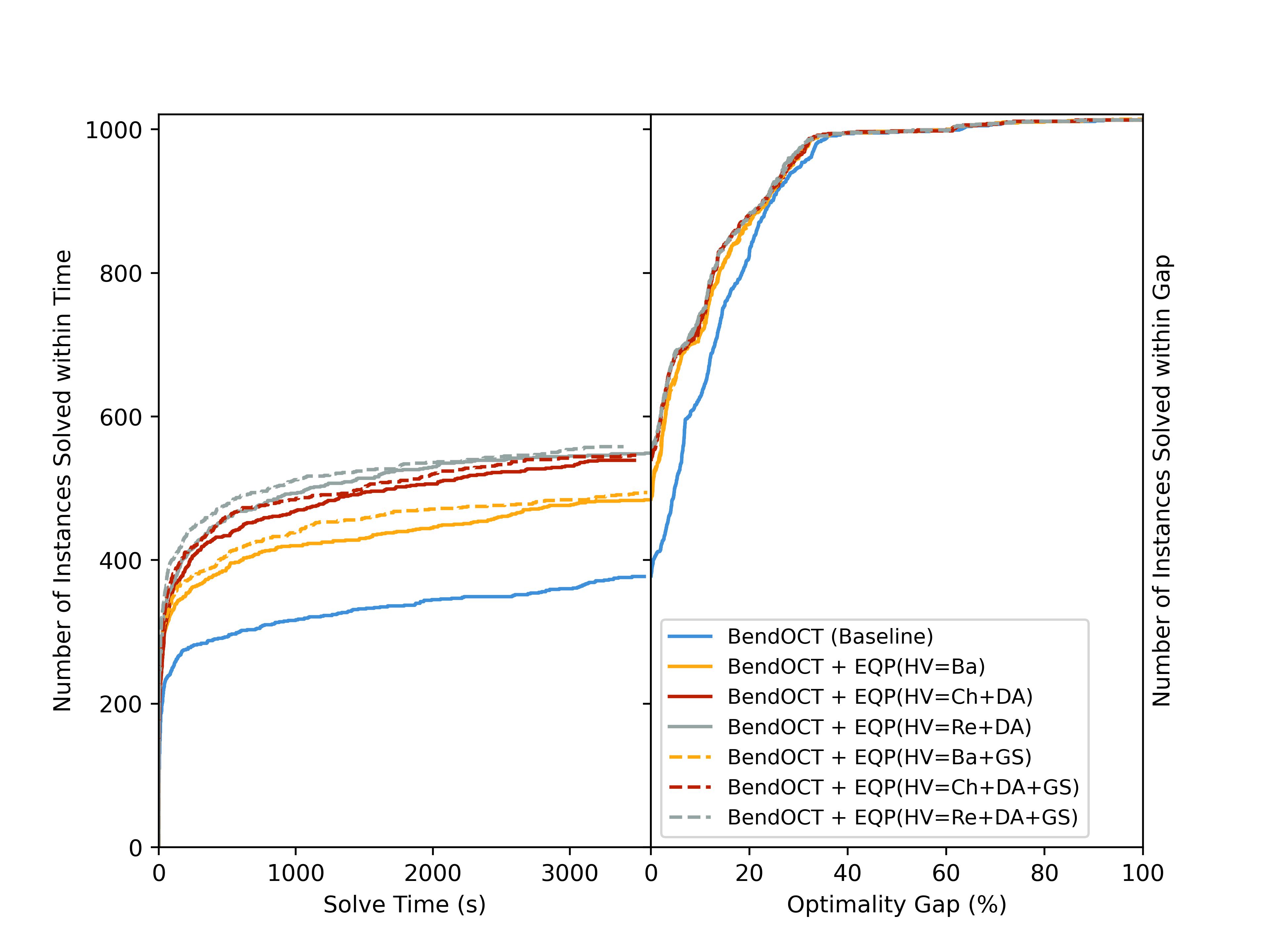}
         \caption{$|F^*| \leq 2$}
    \end{subfigure}
    \caption{Equivalent point inequalities computational results assessing the effect of group selection (GS) constraints. Group selection variants are dashed.}
    \label{fig:EQP Group Selection Comparison}
\end{figure}

Figure \ref{fig:EQP Disaggregation Baseline} assesses the effect of constraint disaggregation on $\mathcal{H}^{Chain}$ and $\mathcal{H}^{Recursive}$. In both cases constraint disaggregation is a notable improvement, and when $|F^*| \leq 2$ we find that $\mathcal{H}^{Recursive}$ maintains an edge over $\mathcal{H}^{Chain}$. Figure \ref{fig:EQP Group Selection Comparison} assesses the effectiveness of group selection constraints $\mathcal{G}^{GS}$ across all variants of $\mathcal {H}$ with constraint disaggregation applied for $\mathcal{H}^{Chain}$ and $\mathcal{H}^{Recursive}$. The results show a marginal improvement across all variants while preserving their relative effectiveness.

\begin{figure}[htb!]
    \centering
    \includegraphics[width=\linewidth]{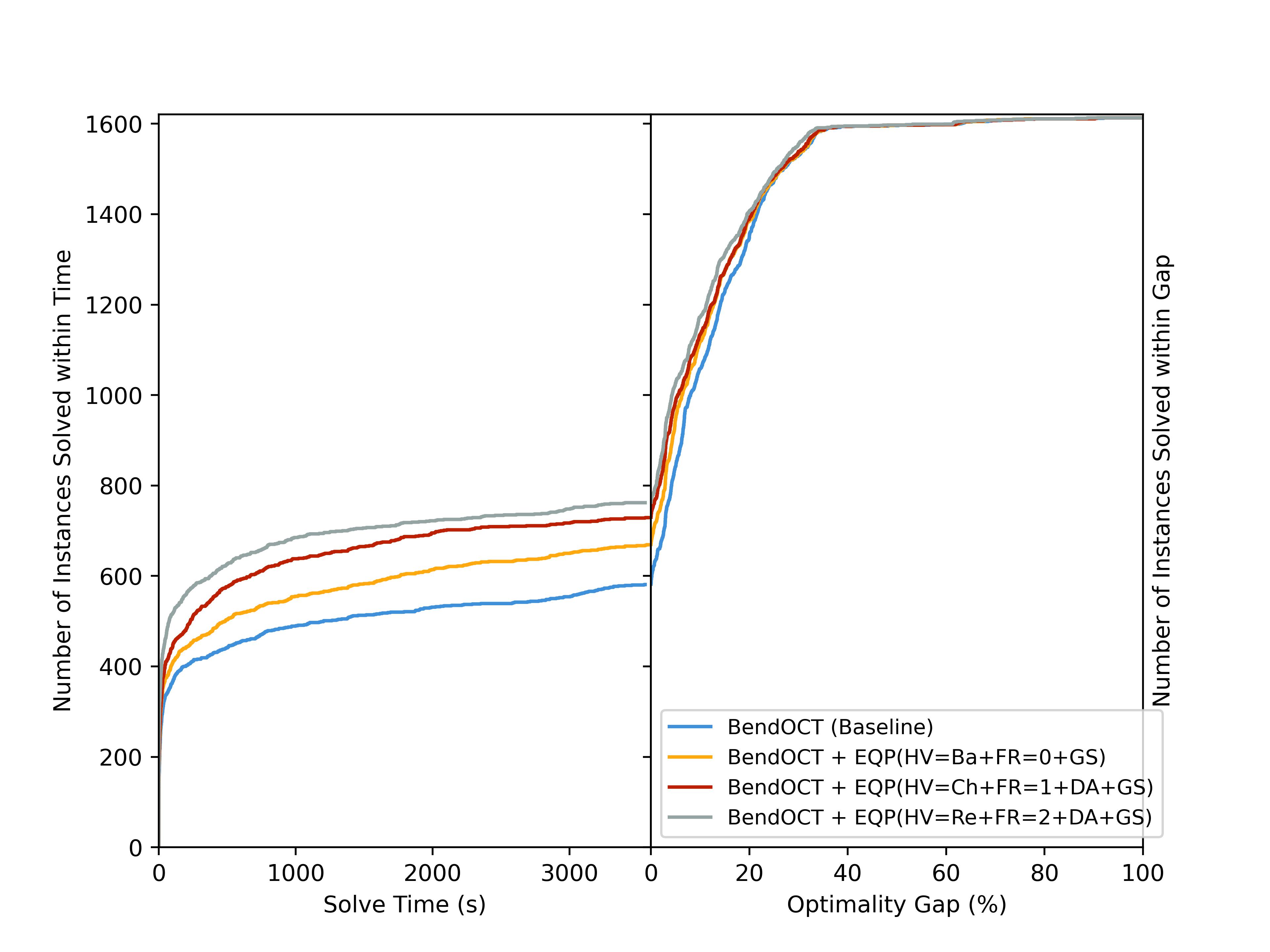}
    \caption{Equivalent point inequalities computational results assessing effectiveness of the best configurations for each upper bound on $|F^*|$. $FR$ indicates the maximum cardinality of $|F^*|$}
    \label{fig:Equivalent Point Best Configurations Comparison}
\end{figure}

While splitting the results according to the maximum cardinality of $F^*$ is convenient for comparing the effect of different configurations, it reveals little about the relative effectiveness of different upper bounds on $|F^*|$. It is not obvious that a higher bound will always yield better performance. As the upper bound increases, all previously modelled EQP sets and associated constraints are included as well as new EQP sets with larger split sets. The number of EQP sets modelled can grow aggressively with the size of the bound for some datasets, and we expect that cuts derived from EQP sets with larger split sets will be less effective since more features can split the samples in the relaxation. Eventually the increased relaxation solve times will outweigh the marginal tightening of the relaxation. Figure \ref{fig:Equivalent Point Best Configurations Comparison} compares the highest performing variants for each cardinality bound. We find that the highest performing configuration is $\mathcal{H}^{Recursive}$ with constraint disaggregation, group selection constraints $\mathcal{G}^{GS}$, and a cardinality bound of $|F^*| \leq 2$.

\subsubsection{Path Bound Cutting Planes}
Section \ref{sect:Path Bound Cutting Planes} describes a number of variants of the basic cutting planes. We present results over combinations of the following:

\begin{itemize}
    \item \textbf{Solution polishing} - We test the effect of the solution polishing primal heuristic as a means of closing the optimality gap when the solver has a tight bound but fails to converge to the optimal solution. We do not check the feasibility of the solutions with respect to the full model before running the primal heuristic.
    \item \textbf{Bounding negative samples (BNS)} - We assess the effect of modifying the cutting planes to force the classification score of samples misclassified in the optimised subtree to zero.
    \item \textbf{Bounding the tree structure (BSt)} - We assess the effect of passing the structure of the optimised subtree to the solver. This is tested in conjunction with both the basic cuts and the BNS variant.
\end{itemize}

\begin{figure}[htb!]
    \centering
    \includegraphics[width=\linewidth]{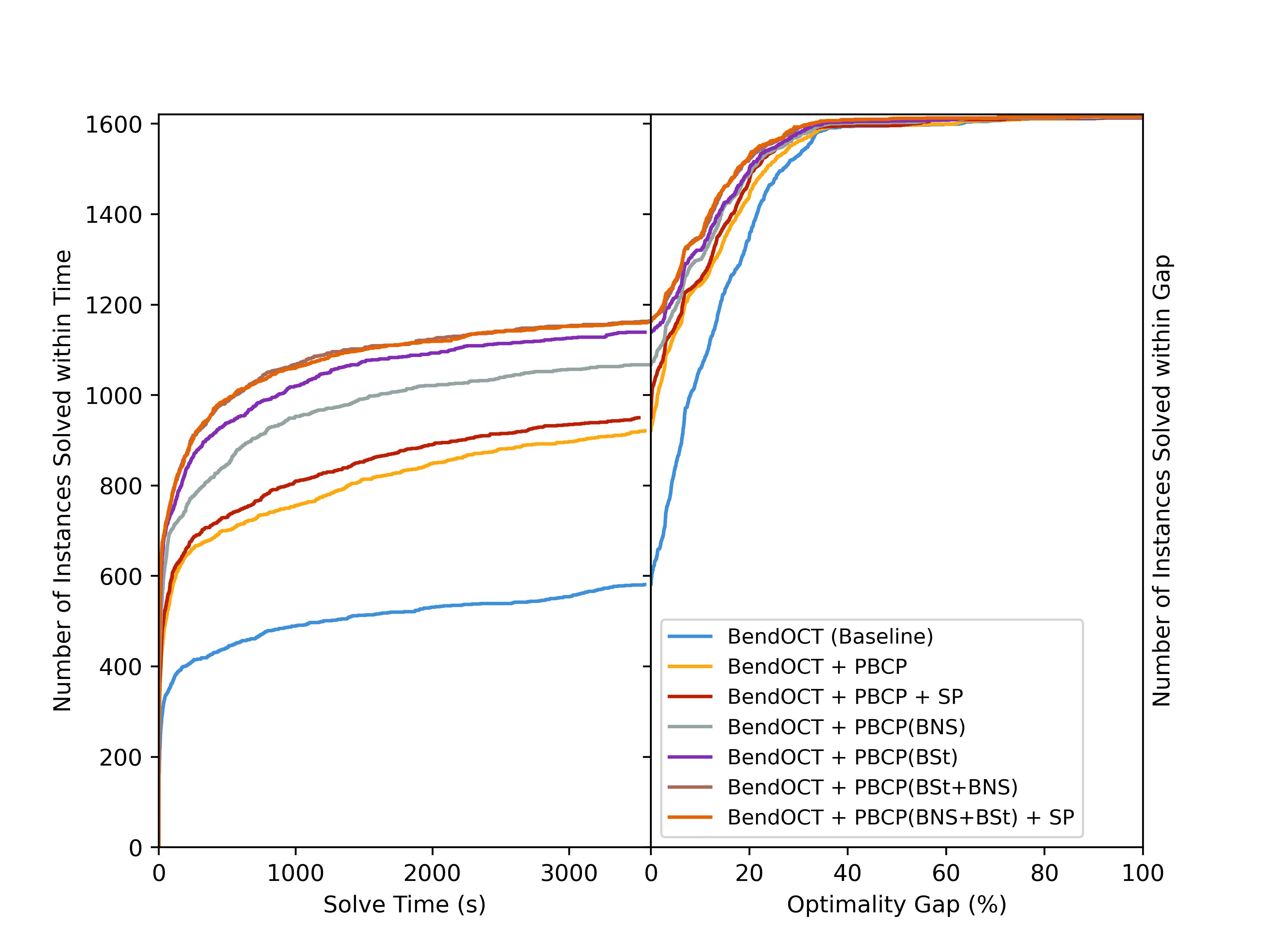}
    \caption{Computational results for variants of path bound cutting planes.}
    \label{fig:PBCP All Variants Comparison}
\end{figure}

Results for all variants with cuts added are shown in Figure \ref{fig:PBCP All Variants Comparison}. The results show that the basic version of the cuts are very effective and that solution polishing, bounding negative samples, and bounding the tree structure each have a positive impact when applied individually with tree structure bounds having the largest impact. The strongest results are found when bounding negative samples and the tree structure. Further applying solution polishing has no effect indicating that the variants are effective at forcing the solver to adopt strong subtree solutions. We report per dataset results for the path bound cutting planes in Table \ref{tab:Aggregated PBCP Results} which reinforce the effectiveness of the cutting planes and the proposed improvements.

\begin{table}[htb]
    \centering
    \renewcommand*{\arraystretch}{0.85}
    \begin{tabular}{cccccccccc}
\hline
\multirow{2}{*}{Dataset} & \multicolumn{3}{c}{BendOCT} & \multicolumn{3}{c}{PBCP} & \multicolumn{3}{c}{PBCP + BSt + BNS}\\
  & Solved & Time & Gap & Solved & Time & Gap & Solved & Time & Gap\\
\hline
soybean-small & $30$ & $0.2$ & - & $30$ & $0.2$ & - & $\mathbf{30}$ & $\mathbf{0.1}$ & \textbf{-}\\
monk3 & $25$ & $258.5$ & $0.30$ & $30$ & $30.7$ & - & $\mathbf{30}$ & $\mathbf{8.3}$ & \textbf{-}\\
monk1 & $30$ & $6.4$ & - & $30$ & $5.1$ & - & $\mathbf{30}$ & $\mathbf{0.7}$ & \textbf{-}\\
hayes-roth & $20$ & $353.6$ & $2.72$ & $30$ & $437.1$ & - & $\mathbf{30}$ & $\mathbf{47.3}$ & \textbf{-}\\
monk2 & $15$ & $1642.6$ & $17.89$ & $27$ & $635.7$ & $1.30$ & $\mathbf{30}$ & $\mathbf{46.1}$ & \textbf{-}\\
house-votes-84 & $26$ & $257.3$ & $0.27$ & $30$ & $31.9$ & - & $\mathbf{30}$ & $\mathbf{10.7}$ & \textbf{-}\\
spect & $18$ & $956.7$ & $3.63$ & $24$ & $359.7$ & $1.15$ & $\mathbf{30}$ & $\mathbf{167.0}$ & \textbf{-}\\
breast-cancer & $7$ & $468.7$ & $16.26$ & $18$ & $299.4$ & $19.93$ & $\mathbf{19}$ & $\mathbf{185.9}$ & $\mathbf{17.38}$\\
balance-scale & $18$ & $785.5$ & $20.54$ & $19$ & $325.7$ & $24.34$ & $\mathbf{22}$ & $\mathbf{454.4}$ & $\mathbf{15.98}$\\
tic-tac-toe & $4$ & $236.7$ & $22.56$ & $13$ & $1413.2$ & $15.11$ & $\mathbf{19}$ & $\mathbf{140.2}$ & $\mathbf{17.01}$\\
car\_evaluation & $6$ & $563.7$ & $14.99$ & $18$ & $123.7$ & $19.15$ & $\mathbf{21}$ & $\mathbf{389.2}$ & $\mathbf{14.25}$\\
kr-vs-kp & $7$ & $1475.1$ & $5.79$ & $17$ & $860.2$ & $5.69$ & $\mathbf{19}$ & $\mathbf{551.3}$ & $\mathbf{5.82}$\\
\hline
iris & $55$ & $243.6$ & $0.73$ & $60$ & $33.0$ & - & $\mathbf{60}$ & $\mathbf{11.3}$ & \textbf{-}\\
wine & $48$ & $464.6$ & $1.17$ & $49$ & $191.0$ & $1.20$ & $\mathbf{51}$ & $\mathbf{168.3}$ & $\mathbf{0.81}$\\
plrx & $11$ & $435.5$ & $22.56$ & $30$ & $615.5$ & $19.73$ & $\mathbf{37}$ & $\mathbf{150.5}$ & $\mathbf{21.81}$\\
wpbc & $9$ & $124.9$ & $16.86$ & $14$ & $348.3$ & $9.70$ & $\mathbf{34}$ & $\mathbf{81.5}$ & $\mathbf{14.76}$\\
parkinsons & $24$ & $351.0$ & $5.74$ & $33$ & $388.5$ & $5.18$ & $\mathbf{36}$ & $\mathbf{47.3}$ & $\mathbf{4.92}$\\
sonar & $3$ & $3055.4$ & $15.64$ & $5$ & $691.8$ & $12.41$ & $\mathbf{22}$ & $\mathbf{562.3}$ & $\mathbf{12.88}$\\
wdbc & $13$ & $298.4$ & $3.95$ & $20$ & $1020.2$ & $3.02$ & $\mathbf{35}$ & $\mathbf{227.0}$ & $\mathbf{3.98}$\\
transfusion & $15$ & $281.3$ & $14.15$ & $38$ & $107.6$ & $16.89$ & $\mathbf{60}$ & $\mathbf{197.0}$ & \textbf{-}\\
banknote & $35$ & $695.9$ & $2.95$ & $44$ & $436.0$ & $3.32$ & $\mathbf{60}$ & $\mathbf{253.7}$ & \textbf{-}\\
ozone-one & $19$ & $295.9$ & $2.37$ & $23$ & $194.0$ & $2.51$ & $\mathbf{41}$ & $\mathbf{885.2}$ & $\mathbf{2.54}$\\
segmentation & $0$ & - & $50.47$ & $0$ & - & $47.44$ & $\mathbf{25}$ & $\mathbf{1555.5}$ & $\mathbf{41.44}$\\
spambase & $0$ & - & $15.60$ & $5$ & $1789.0$ & $14.62$ & $\mathbf{20}$ & $\mathbf{2120.5}$ & $\mathbf{13.27}$\\
\hline
hepatitis & $33$ & $1145.8$ & $0.78$ & $59$ & $558.6$ & $0.25$ & $\mathbf{60}$ & $\mathbf{246.0}$ & \textbf{-}\\
fertility & $31$ & $619.6$ & $2.88$ & $60$ & $696.3$ & - & $\mathbf{60}$ & $\mathbf{188.0}$ & \textbf{-}\\
ionosphere & $10$ & $584.6$ & $8.36$ & $19$ & $1003.3$ & $5.67$ & $\mathbf{35}$ & $\mathbf{124.3}$ & $\mathbf{7.24}$\\
thoracic & $14$ & $156.3$ & $11.69$ & $35$ & $605.6$ & $12.35$ & $\mathbf{40}$ & $\mathbf{211.4}$ & $\mathbf{11.92}$\\
ILPD & $10$ & $419.8$ & $29.31$ & $29$ & $690.7$ & $25.58$ & $\mathbf{36}$ & $\mathbf{109.0}$ & $\mathbf{29.49}$\\
credit & $12$ & $90.0$ & $12.09$ & $33$ & $1120.5$ & $11.77$ & $\mathbf{36}$ & $\mathbf{167.7}$ & $\mathbf{11.81}$\\
biodeg & $0$ & - & $24.54$ & $7$ & $870.7$ & $21.48$ & $\mathbf{30}$ & $\mathbf{624.5}$ & $\mathbf{20.94}$\\
seismic-bumps & $16$ & $9.1$ & $5.55$ & $37$ & $154.5$ & $5.24$ & $\mathbf{39}$ & $\mathbf{67.2}$ & $\mathbf{5.58}$\\
ann-thyroid & $18$ & $357.4$ & $5.66$ & $36$ & $927.2$ & $5.59$ & $\mathbf{38}$ & $\mathbf{229.1}$ & $\mathbf{5.94}$\\
\hline
\end{tabular}
    \caption{Path bound cutting planes per dataset results. The solved columns list the number of instances solved to optimality. Time lists the averages solve time in seconds of solved instances, ’-’ indicates that no instances were solved. Gap lists the average percentage optimality gap in unsolved instances, ’-’ indicates that all instances were solved. Bold font indicates the best performing variant prioritising the number of instances solve to optimality, the average solve time, followed by the average optimality gap.}
    \label{tab:Aggregated PBCP Results}
\end{table}

\end{document}